%% file: minor_revision.tex
\documentclass[hidelinks,onefignum,onetabnum]{siamart220329}

\input{ex_shared}

\ifpdf
\hypersetup{
  pdftitle={Convergence rates of SOS hierarchies for polynomial semidefinite programs},
  pdfauthor={Hoang Anh Tran, \ 
  Kim-Chuan Toh}
}
\fi

\DeclareMathOperator{\lb}{lb}
\DeclareMathOperator{\dist}{dist}

\def\inprod#1#2{\langle#1,\,#2\rangle}

\allowdisplaybreaks
\begin{document}

\maketitle

\begin{abstract}
We introduce a novel moment-SOS hierarchy of lower bounds for a polynomial optimization problem (POP) whose feasible set is defined by polynomial matrix inequalities (PMIs). Our hierarchy avoids the Kronecker product structure in Hol-Scherer's hierarchy, thus resulting in smaller-sized semidefinite programs. Our approach involves utilizing a penalty function framework to directly address the matrix-based constraint, which is applicable to both discrete and continuous polynomial optimization problems. We investigate the convergence rates of these bounds for both types of problems. The proposed method yields a variant of Putinar's theorem, tailored for positive polynomials on a compact set $\mathcal{X}$ defined by a polynomial matrix inequality. More specifically, we derive novel insights into the bounds on the degree of the SOS polynomials required to certify positivity over $\mathcal{X}$, based on Jackson's theorem and a variant of the Łojasiewicz inequality in the matrix setting.

\end{abstract}

% REQUIRED

% REQUIRED

\section{Introduction}
This paper considers solving a polynomial optimization problem (POP) with a polynomial matrix semidefinite inequality of the following form: 
\begin{align} 
    f_{\min} &= \min{f(x)}\nonumber\\
	& \mbox{subject to } x \in \mathcal{X} := \{ x\in \mathbb{R}^n \,:\,
	 G(x)\succeq 0\}, \label{SPOP}
\end{align}
where $f \in \mathbb{R}[x]$ is a polynomial of degree $d$, and $G(x)=(g_{ij}(x))_{m \times m}$ is an $m \times m$ symmetric polynomial matrix with $g_{ij}(x)\in  \mathbb{R}[x]$. 
In this paper, 
we use $\mathbb{R}[x]$ to denote the ring of polynomials in $x$, and
$\bS\bR[x]^m$ to denote the set of $m \times m$  symmetric polynomial  matrices. 
The  set of polynomials of degree at most $k$ is denoted by $\bR[x]_k$. 
The degree of $G(x)$, denoted as $\mathrm{deg}(G)$,  
is defined as the maximum degree among all its entries $g_{ij}(x)$.
Throughout this paper, we denote the degree of $f$ by $d$, and 
set $l = \lceil \deg(G)/2 \rceil$.

As a special case of \eqref{SPOP}, a classical polynomial optimization problem corresponds to $G(x)$ being a diagonal matrix with polynomials $g_{ii}(x)$, $i\in[n]$, on its diagonal. 
 In this case, $\mathcal{X}$ is a semialgebraic set defined by polynomial inequalities $g_{ii}(x) \geq 0$, $i\in[n]$. Problem $\eqref{SPOP}$ extends beyond this special
 case by specifying the feasible set through a polynomial matrix inequality. Throughout this paper, we refer to \eqref{SPOP} as a matrix POP to distinguish it from the scalar POP, where the feasible set is defined by polynomial inequalities.

One may ask whether it is possible to reformulate $\mathcal{X}$ as a semialgebraic set.
By denoting the eigenvalues of $G(x)$ for $x\in\mathbb{R}^n$ by $\lambda_1(x) \geq \dots \geq \lambda_m(x)$, we get the following equivalent description of the feasible set of \eqref{SPOP}:
\begin{align}\label{semialg by semid}
	\mathcal{X} =\;   \{x \in \mathbb{R}^n: \lambda_m(x) \geq 0\}.
\end{align} 
However, we should note that since the function $\lambda_m(x)$ is generally not a polynomial, it cannot represent $\mathcal{X}$ as a semialgebraic set. 
To represent $\mathcal{X}$ as a semialgebraic set with polynomial inequalities, a standard method is to apply Descartes' rule of signs, which explicitly describes $\mathcal{X}$ via $m$ scalar polynomial inequalities involving the coefficient polynomials $c_i(x)$ of 
the characteristic polynomial of 
$G(x)$  defined as follows (see, e.g.,\cite{henrion2006convergent,lasserre2005unified}):
$
\det(tI_m - G(x))=\; t^m +\sum_{i=1}^m (-1)^ic_i(x)t^{m-i}, \; t\in \mathbb{R}.
$
Based on a generalization of Descartes' rule of signs (see, e.g., \cite{lasserre2005unified}), the polynomial matrix inequality $G(x) \succeq 0$ is equivalent to a system of polynomial inequalities $c_i(x) \geq 0, \; \forall i \in [m]$. That is, $\mathcal{X}$ admits the form of a basic semialgebraic set defined as 
$\mathcal{X} =\; \{x \in \mathbb{R}^n: \ c_i(x) \geq 0\ \forall i = 1,\dots, m \}.$
However, the degree of $c_i$ can be significantly larger than the degree of $G$. 
As an illustration, consider the case where all entries $g_{ij}$ have the same degree $\ell$,
and hence the degree of $G$ is $\ell$, but the degree of $c_1$ is $(m-1)\ell$, which is much larger than $\ell$. Furthermore, computing all $c_i$'s is  computationally as expensive as calculating the determinant of $tI_m- G(x)$, which require up to $m!$ steps. Consequently, reformulating \eqref{SPOP} as a classical POP and applying the corresponding moment-SOS hierarchy can be prohibitively expensive. Additionally, since the definition of a quadratic module depends on the description of the underlying semialgebraic set, distinct quadratic modules are obtained depending on whether the underlying set is defined by polynomial inequalities or a matrix inequality. Thus, constructing a hierarchy that operates directly on the matrix inequality is computationally advantageous and practically preferable.

\subsection{Hierarchy of lower bounds in the matrix case}
\label{sec-1.2}

We first review the works of Hol and Scherer \cite{hol2004sum, hol2005sum} and Kojima \cite{Kojima} who proposed a moment-SOS hierarchy based on SOS polynomial matrix  for \eqref{SPOP}.

A polynomial $p \in \mathbb{R}[x]$ is said to be a {\em sum-of-squares (SOS)} if it can be written as $p =p_1^2 +\ldots +p_k^2$ for some $p_i \in \mathbb{R}[x]$, $i\in[k]$. 
Extending this to the matrix case, we say that a symmetric $m \times m$ polynomial matrix $P(x)$ is a {\em SOS polynomial matrix} if there exists a (not necessarily square) polynomial matrix $T(x)$ such that $P(x) = T(x)^{\top}T(x)$. We denote the set of all SOS polynomials and the set of all $m \times m$ SOS polynomial matrix by $\Sigma[x]$ and $\bS\Sigma[x]^m$, respectively. We define the {\em matrix quadratic module} $\cQ(\cX)$ associated with $\cX$ in \eqref{SPOP} as: 
\begin{equation} \label{eq-qmodule}
    \cQ(\mathcal{X}):= \left\{ \sigma(x) + \langle R(x),G(x)\rangle \;:\; \sigma \in \Sigma[x],\ R(x) \in \bS\Sigma[x] \right\}.
\end{equation}
In the above, $\inprod{\cdot}{\cdot}$ denotes the trace inner product between two symmetric matrices. 
It is worth noting that $\cQ(\cX)$ depends on the description 
\eqref{SPOP} of $\cX$ rather than the set $\cX$ itself. Under the Archimedean condition, Hol and Scherer established a natural extension of Putinar's Positivstellensatz \cite{putinar1993positive}, as stated in the following theorem.

\begin{theorem}\cite[Theorem 1]{hol2004sum}\label{GPutinar} 
Suppose $\mathcal{X}$ satisfies the Archimedean condition, that is, there exist an SOS polynomial $\sigma$, an SOS polynomial matrix $R(x)$, and a scalar $N$ such that 
\begin{displaymath}
    N - x^\top x  = \sigma(x) +  \langle R(x),G(x) \rangle.
\end{displaymath}
Then every positive polynomial $f$ on $\mathcal{X}$ belongs to the quadratic module $Q(\mathcal{X})$.
\end{theorem}

Next, we consider the {\em truncated} quadratic module, defined as follows: 
\begin{multline*}
    \cQ(\cX)_{2r}:= \Bigl\{ \sigma(x) + \langle R(x),G(x)\rangle\;: \sigma \in \Sigma[x],\ R(x) \in \bS\Sigma[x],\\ \deg \sigma \leq 2r,\; \deg R + \deg G \leq 2r \Bigr\},
\end{multline*}
whose membership can be checked by a semidefinite program (SDP) (see e.g., \cite{hol2004sum,hol2005sum,henrion2006convergent}). Whence, Hol and Scherer proposed the following hierarchy of lower bounds to approximate the optimal value $f_{\min}$ of $f$ over $\cX$:
\begin{equation}\label{hierarchy: Hol and Scherer}
    \lb(f,\cQ(\cX))_r :=\; \max\{t \in \bR\;:\; f- t \in \cQ(\cX)_{2r}\}.
\end{equation}
Under the Archimedean condition on $\cX$, the lower bound $\lb(f,\cQ(\cX))_r$ converges to $f_{\min}$ as $r \to \infty$. 

Computing $\lb(f,\cQ(\cX))_r$ can be cast as an SDP. In the scalar case, the primary disadvantage is that the SDP at the $r$-th level of the moment-SOS hierarchy has $\binom{n+2r}{n}$ variables, with the size of the SDP matrix variable equals to $\binom{n+r}{n}$. This limitation becomes even more severe for 
the matrix moment-SOS hierarchy, where the size of the SDP matrix variable is $m \cdot\binom{n+r}{n}$, due to the Kronecker product structure involved in the hierarchy (this is elaborated in Remark~\ref{rem: exploxion of degree}).

 %%%%%%%%%%%%%%%%%%%%%%%%%%%
\subsection{Related works}
    Polynomial optimization with polynomial matrix inequality has  applications in various fields such as control theory \cite{chesi2010lmi,henrion2006convergent,ichihara2009optimal, scherer2006lmi}, quantum information theory \cite{fang2021sum}, and statistics \cite{henrion2020moment}. Various matrix moment-SOS hierarchies have been introduced and studied in the works \cite{hol2004sum,hol2005sum,henrion2006convergent,Kojima} to address the problem \eqref{SPOP}.

    For the scalar POP, Lasserre \cite{lasserre2001global, lasserre2009moments,lasserre2011new,lasserre2001explicit,lasserre2005unified} introduced a moment-SOS hierarchy and proved its asymptotic convergence
    to the optimal value. He applied this framework to POP over the binary hypercube, which can be described by polynomial equalities (see e.g., \cite{lasserre2016max, lasserre2001explicit}). The finite convergence of the scalar moment-SOS hierarchy was proved under convexity \cite{de2011lasserre}, and some other regular assumptions (see e.g., \cite{huang2023optimality,huang2024finite, marshall2006representations,nie2014optimality}).

    On the convergence rate of the scalar moment-SOS hierarchy under the Archimedean condition, Nie and Schweighofer \cite{nie2007complexity} proved an exponential rate when the order of the hierarchy increases to infinity. Recently, Baldi and Mourrain \cite{baldi2021moment,baldi2025lojasiewicz} improved the exponent rate to polynomial rate. Huang \cite{huang2025complexity} extended the convergence rate in the scalar case \cite{baldi2025lojasiewicz} to the matrix case. As we shall  see later, independent of Huang’s work, our work derives a polynomial rate for a simplified matrix hierarchy. 

    When the feasible set of a POP has special structure, stronger convergence rates are known. The moment-SOS hierarchy for POP over the binary hypercube was shown to have the convergence rate of $\mathcal{O}(1/r^2)$ in \cite{slot2023bpop}. The work \cite{slot2111sum} studied the convergence rates of $\mathcal{O}(1/r^2)$ for the $n$-dimensional unit ball and the standard simplex. The same rate was also proved for the hyper-sphere with homogeneous polynomial objective function in  \cite{fang2021sum}. For the sake of readability, we postpone stating these convergence rates precisely until we have defined the necessary notation 
    in Section~\ref{section: preliminary}.

    For the matrix case, other than Huang's work, we are not aware
    of other works that study the convergence rates of matrix
    moment-SOS hierarchies for a matrix POP.

\subsection{Contribution} 
    This work addresses two fundamental challenges when solving \eqref{SPOP}: (i) We develop a streamlined version of the matrix moment-SOS hierarchy \eqref{hierarchy: Hol and Scherer} that avoids using the Kronecker product structure and hence the size of the SDP matrix variable is not 
    increased by a factor of  $m$; (ii) We rigorously analyze the convergence rate of 
 our proposed hierarchy through a novel penalty function framework.

    In what follows, we define our hierarchy and highlight the main results of this paper. Throughout the paper, for any univariate polynomial $h(t)= h_kt^k+\ldots+h_1t+h_0 \in \bR[t]$ and $G(x) \in \bS\bR[x]^m$, we define the polynomial matrix $h(G(x))$ by 
    \begin{equation}\label{def: univariate poly with matrix}
        h(G(x))\; = \; h_k (G(x))^k+\ldots+h_1G(x)+h_0I_m.
    \end{equation}
   We define the following sets depending on $G$ by
    \begin{align*}
        \cH(G)&=\; \{\langle h(G(x)),G(x) \rangle:\; h(t) \geq 0 \; \forall t \in [-1,1],\; h\in\bR[t]\},\\
        \cH(G)_k&=\; \{\langle h(G(x)),G(x) \rangle:\; h(t) \geq 0 \; \forall t \in [-1,1],\; h \in \bR[t]_k\}.
    \end{align*}
 We propose the following sets as a replacement of the quadratic module in \eqref{hierarchy: Hol and Scherer} for two cases where $\cX$ is contained in the discrete binary hypercube $\bB^n$ and the $n$-dimensional unit ball $B^n$, respectively. 
    \begin{enumerate}
        \item For $\cX \subset \bB^n$, define for $r\in \mathbb{N}$ and $r\geq l$,
            \begin{eqnarray}
                \widetilde{\cQ}(\cX) \;=\; \cQ(\bB^n)+\cH(G), \quad \widetilde{\cQ}(\cX)_{2r} \;=\; \cQ(\bB^n)_{2r}+\cH(G)_{\lfloor r/l \rfloor-1}.
                \label{bpop-qmodule}
            \end{eqnarray}
            The corresponding hierarchy for \eqref{SPOP} is 
        \begin{equation}\label{proposed hierarchy for bpop}
            \lb(f, \widetilde{\cQ}(\cX))_r \;:=\; \max\{t \in \bR:\; 
            f-t \in  \widetilde{\cQ}(\cX)_{2r}\}.
        \end{equation}
        
        \item For $\cX \subset B^n$, define for $r\in \mathbb{N}$ and $r\geq l$,
            \begin{displaymath}
                \overline{\cQ}(\cX) \;=\; \cQ(B^n)+\cH(G), \quad \overline{\cQ}(\cX)_{2r} \;=\; \cQ(B^n)_{2r}+\cH(G)_{\lfloor r/l \rfloor-1}.
            \end{displaymath}
       The corresponding hierarchy for \eqref{SPOP} is 
        \begin{equation}\label{proposed hierarchy for cpop}
            \lb(f, \overline{\cQ}(\cX))_r \;:=\; \max\{t \in \bR:\; f-t \in  \overline{\cQ}(\cX)_{2r}\}.
        \end{equation}
    \end{enumerate}
 Note that the subscript $\lfloor r/l \rfloor-1$ in $\cH(G)_{\lfloor r/l \rfloor-1}$ specifies 
 the bound on the degree of its underlying polynomials $h$, and hence
  $\deg ( \langle h(G(x)),G(x) \rangle) \leq 2l \cdot (\lfloor r/l \rfloor-1)+2l \leq 2r$.

Our contributions consist of the following results: 
\begin{enumerate}
    \item When $\cX$ is a subset of the binary hypercube $\bB^n$, Corollary~\ref{cor: proposed momemnt SDP for bpop} describes a novel SDP relaxation via moment matrices for \eqref{SPOP}, which also points out how our hierarchy reduces the size of SDP matrices compared to those in the Hol-Scherer hierarchy. Theorem~\ref{thm: convergence rate for bpop} provides the tightness of our proposed hierarchy of lower bounds in \eqref{proposed hierarchy for bpop}.
    \item When $\cX$ is a compact subset  of the unit ball, the novel SDP relaxation via moment matrices and the convergence rate of the 
    hierarchy \eqref{proposed hierarchy for cpop} are presented in Corollary~\ref{cor: proposed momemnt SDP for cpop} and Theorem~\ref{convergence_thm_cpop}, respectively. 
    \item Theorem~\ref{thm: bound on degree} states a new Positivstellensatz for positive polynomials over $\cX$ that is 
    contained in $B^n$ and is defined by a matrix inequality as in \eqref{SPOP}. 
\end{enumerate}

We organize the paper according to our results. Section~\ref{section: preliminary} reviews some preliminary concepts and define all notation we use in the paper. Section~\ref{section: BPOP} and Section~\ref{section: Continuous} study the convergence rates of the hierarchies 
\eqref{proposed hierarchy for bpop} and
\eqref{proposed hierarchy for cpop}
for the discrete setting and continuous setting, respectively. 
Section~\ref{section: another version} introduces and proves a new Positivstellensatz for positive polynomials over $\cX$ defined by a matrix inequality as in \eqref{SPOP}. Some technical proofs are postponed to the Appendices.   

\section{Preliminary}\label{section: preliminary}

We first define some notation which will be used in this paper.   
Denote by $\bS^m$ and $\bS^m_+$ the set of $m \times m$ symmetric matrices and the set of $m \times m$ positive semidefinite matrices, respectively.
For any vector $x =(x_1,\dots,x_n) \in \mathbb{R}^n$ and multi-index $\alpha=(\alpha_1,\dots,\alpha_n) \in \mathbb{N}^n$, the monomial $x_1^{\alpha_1}\cdots x_m^{\alpha_n}$ is denoted by $x^{\alpha}$. Consequently, any polynomial $g(x)$ can be written in the form $g(x) = \sum_{\alpha \in \mathbb{N}^n}g_{\alpha}x^{\alpha}$. 
To quantify the size of a polynomial, we define the norm $\|g\|_{\infty} = \sum_{\alpha \in \mathbb{N}^n}|g_{\alpha}|$. 
We use $\bN^n_r$ to denote the set of multi-indices of length at most $r$ with its cardinality 
 equals to $s(n,r) := \binom{n+r}{n}$. For any matrix $M \in \bR^{m \times m}$, $\tr(N)$ denotes its trace. For two matrices $M_1$ and $M_2$ of the same size, $\langle M_1,M_2 \rangle$ denotes the inner product $\tr(M_1^{\top}M_2)$. 
 For any $t \in \bR$, we denote by $\lceil t \rceil$ (resp. $\lfloor t  \rfloor$) the smallest (resp. largest) integer that is no smaller (resp. larger) than $t$. 
 
 For any matrices $M$ and $N$, we denote their Kronecker product by $M \otimes N$. For any set $S$ and a real-valued function $f$ defined on it, we define $\|f\|_S= \max_{x \in S}|f(x)|$. For any matrix $M \in \bS^m$, we denote its spectral radius by $\rho(M)$. For a set $S \subset \bR^n$, we denote the set of non-negative polynomials over $S$ by 
 $\cP_+(S)$ and those
 of degree at most $k$ by $\cP_+(S)_k$. We use the notation $\dist(\cdot,\cdot)$ to denote the Euclidean distance between $2$ points or between a point and a set in an Euclidean space. 

%%%%%%%%%%%%%%%%%%%%
\subsection{Hierarchy of lower bounds in the scalar case}

We review Lasserre's hierarchy for a scalar POP over a basic semialgebraic set. To distinguish between different constraint types, we use the notation $\mathcal{X}$ to denote a set defined by a matrix inequality, and $\mathit{X}$ to denote a set defined by polynomial inequalities. Let $\mathit{X} \in \mathbb{R}^n$ be a semialgebraic set of the form
\begin{align}\label{semialg by poly}
	\mathit{X} := \{x \in \mathbb{R}^n: g_j(x) \geq 0 \ \forall j \in [m]\},
\end{align} 
where $g_j \in \mathbb{R}[x]$ is a polynomial for all $j \in [m]$. 
Consider the problem
\begin{align}\label{pop}
	f_{\min}^X := \min_{x \in \mathit{X}} f(x)= \max\{ t \in \mathbb{R}: f- t \in  \mathcal{P}_+(\mathit{X})\}.
\end{align}
The quadratic module corresponding to $\mathit{X}$ is defined as
\begin{displaymath}
	\mathcal{Q}(\mathit{X}) :=\; \left\{\sum_{i=0}^{m}\sigma_ig_i: \sigma_i \in \Sigma[x]\right\} \quad ( \text{ where } g_{0}:=1).
\end{displaymath}
For $r \in \bN$, we define the set of SOS polynomials of degree at most $2r$ and the truncated quadratic module as: 
\begin{displaymath}
    \Sigma[x]_{2r} :=\;\Sigma[x] \cap \mathbb{R}[x]_{2r},\quad
                    \mathcal{Q}(\mathit{X})_{2r} :=\; \left\{\sum_{i=0}^{m}\sigma_ig_i: \sigma_i \in \Sigma[x]_{2r},
                    \ \deg \sigma_ig_i \leq 2r\right\}.
\end{displaymath}
It is clear from the definitions that
\begin{displaymath}
	    \Sigma[x] \subseteq\; \mathcal{Q}(\mathit{X}) \subseteq\; \mathcal{P}_+(\mathit{X}), \quad \Sigma[x]_{2r} \subseteq\; \mathcal{Q}(\mathit{X})_{2r} \subseteq\;  \mathcal{P}_+(\mathit{X}).
\end{displaymath}
    Therefore,  approximating $\mathcal{P}_+(X)$ by $\mathcal{Q}(\mathit{X})_{2r}$
    leads to the following scalar hierarchy of lower bounds for the optimal value $f_{\min}^X$ of \eqref{pop}:
	\begin{equation}\label{scalar hierarchy}
	    \lb(f,\mathcal{Q}(\mathit{X}))_r :=\; \max\{t \in \mathbb{R}: f- t \in \mathcal{Q}(\mathit{X})_{2r}\}.
	\end{equation}
    This hierarchy is commonly referred to as the Putinar-type hierarchy of lower bounds. While another hierarchy—the Schmüdgen-type hierarchy—exists (see, e.g., \cite{lasserre2009moments}), we focus exclusively on the former for simplicity. We introduce the following feasible sets, for which we will utilize results of their corresponding Putinar-type hierarchies of lower bounds in our proofs later:
    \begin{align*}
        & \bB^n \;:=\; \{0,1\}^n\;=\; \{ x \in \bR^n\;:\; x_i^2-x_i=0\},\quad 
         B^n \;:=\; \{x \in \bR^n\;:\; \sum_{i=1}^n x_i^2 \leq 1\},\\
        & [0,1]^n \;:=\; \{ x \in \bR^n\;:\; 0 \leq x_i \leq 1 \; \forall i \in [n]\}.
    \end{align*}
    
    %%%%%%%%%%%%%%%%
    \subsection{Moment matrix and localizing matrix}
	We define the moment matrix and localizing matrix for both types of semialgebraic sets: those defined by  polynomial inequalities as in \eqref{semialg by poly} and those defined by a polynomial matrix inequality as in \eqref{SPOP}. Let $b(x)$ denote the canonical basis 
	of monomials in $x$, given by
    \begin{equation}\label{monomial basis}
        b(x) = (1,\; x_1,\ldots,\; x_n,\; x_1^2,\; x_1x_2,\;\ldots,\;x_1x_n,\;\ldots)^{\top}= (x^{\alpha})^{\top}_{\alpha \in \bN^n}.
    \end{equation}
    For any $r \in \bN$, we denote the canonical basis of monomials with degree at most 
    $r$ by $b_r(x) = (x^{\alpha})^{\top}_{\alpha \in \bN^n_r}$. For a fixed sequence $y= (y_{\alpha})_{\alpha \in \bN^n}$ indexed by the basis $b(x)$, we define the Riesz functional $L_{y}\;:\; \bR[x] \to \bR$ by 
    \begin{displaymath}
        f(x) \;=\; \sum_{\alpha \in \bN^n}f_{\alpha}x^{\alpha} \quad \mapsto \quad L_y(f)\;=\; \sum_{\alpha \in \bN^n}f_{\alpha}y_{\alpha}.
    \end{displaymath}
	The {\em moment matrix} $\M(y)$ is an infinite matrix with rows and columns indexed by the basis $b(x)$. For $\alpha, \beta \in \bN^n$, the $(\alpha,\beta)$-entry  of $\M(y)$ is defined by 
    \begin{displaymath}
        [\M(y)]_{\alpha,\beta}\;=\; L_y([b(x)b(x)^{\top}]_{\alpha,\beta})\;=\; 
        L_y(x^{\alpha+\beta}) \;=\; y_{\alpha+\beta}.
    \end{displaymath}
    We denote by $\M_r(y)$ the finite {\em truncation} of $\M(y)$ with rows and columns indexed by the basis $b_r(x)$. 
    \begin{comment}
    For a polynomial $g = \sum_{\alpha \in \bN^n}g_{\alpha}x^{\alpha}$. We denote by $\M_r(gy)$ the {\em localizing} matrix associated with $y$ and $g$, defined by 
    \begin{displaymath}
        [\M_r(gy)]_{\alpha\beta}\;=\; L_y(g(x)[b(x)b(x)^{\top}]_{\alpha,\beta})\;=\; \sum_{\gamma \in \bN^n}g_{\gamma}y_{\alpha+\beta+\gamma} \quad \forall \alpha,\; \beta \in \bN^n_r.
    \end{displaymath}
    In the description \eqref{semialg by poly} of $X$, we set $d_j := \lceil \deg(g_j)/2 \rceil$  for $j \in [m]$. Then for any $2r \geq \max\{d, 2d_1,\ldots,2d_m\}$, consider the following SDP relaxation of \eqref{pop}:
    \begin{align}\label{primal SDP for pop}
        \min_{y \in \bR^{s(n,2r)}} & \quad L_y(f)\\
        \mbox{subject to}& \quad y_0 =1,\; \M_r(y) \succeq 0,\; \M_{r-d_j}(g_jy) \succeq 0 \quad \forall j \in [m].\nonumber
    \end{align}
    This SDP relaxation is also referred to the SDP relaxation via {\em moment} of \eqref{pop}. The pair \eqref{primal SDP for pop} and \eqref{scalar hierarchy} forms a primal-dual semidefinite programming, for which the strong duality holds under the Archimedean condition of $\mathit{X}$ (see e.g., \cite{josz2016strong}). We notice that \eqref{primal SDP for pop} has $s(n,2r)$ variables the SDP matrix of the size $s(n,r)$ corresponding to the size of the truncation moment matrix $\M_r(y)$. 
    \end{comment}
    
    Henrion and Lasserre \cite{henrion2006convergent} studied the SDP relaxation of \eqref{SPOP} via moment relaxation, which combines with \eqref{hierarchy: Hol and Scherer} to form a primal-dual SDP. Their approach extends the definition of a localizing matrix to the polynomial matrix case via taking Kronecker product. For the description  $\cX$ in \eqref{SPOP}, we define the localizing matrix $\M_r(Gy)$ by
    \begin{displaymath}
        \M_r(Gy) \;=\; L_y\big(G(x) \otimes (b_r(x)b_r(x)^{\top})\big),
    \end{displaymath}
    where we slightly abuse the notation of the Riesz functional to mean that $L_y$ acts entry-wise on the polynomial matrix $G(x) \otimes (b_r(x)b_r(x)^{\top})$. Therefore, for any $2r \geq \max\{d,\;l\}$, the  moment relaxation of \eqref{SPOP} is given by 
    \begin{align}\label{primal SDP for SPOP}
        \min_{y \in \bR^{s(n,2r)}} & \Big\{ L_y(f)\,:\,
        %\mbox{subject to}& \quad 
        y_0 =1,\; \M_r(y) \succeq 0,\; \M_{r-l}(Gy) \succeq 0\Big\}.
        \nonumber
    \end{align}
    When $G$ is a scalar polynomial $g$, i.e., $m=1$, we have $\M_{r-l}(Gy)=\M_{r-l}(gy)$, which is the standard definition of the localizaing matrix  for $g$, and \eqref{primal SDP for SPOP} becomes the moment relaxation for a scalar POP (we refer to \cite{lasserre2009moments} for the standard notation of Lasserre's hierarchy for a scalar POP).
    
    \begin{remark}\label{rem: exploxion of degree}
        While the number of variables for $y$ 
        in \eqref{primal SDP for SPOP} remains at $\binom{n+2r}{2r}$, the SDP matrix size grows to $m\cdot \binom{n+r-l}{n}$, which becomes prohibitively large compared to the SDP matrix size of $\binom{n+r}{r}$ in the scalar case. Hol and Scherer \cite{hol2004sum} proposed an alternative SDP relaxation, but it still employs the Kronecker product operation, thus leading to a similar explosion in the dimension. This limitation motivates our work to develop a novel variant of \eqref{hierarchy: Hol and Scherer} or \eqref{primal SDP for SPOP}, that is either \eqref{proposed hierarchy for bpop} or \eqref{proposed hierarchy for cpop}, which avoids the Kronecker product structure, thereby preventing the rapid growth in the SDP matrix size.
    \end{remark}

    \subsection{Convergence rates of moment-SOS hierarchies}\label{sec: summary on complexity}
    %We summarize the convergence rate results for the moment-SOS hierarchy in both scalar and matrix cases. 
    
    We first list  the convergence rates of the moment-SOS hierarchies for problem \eqref{pop} over the binary 
    hypercube $\mathbb{B}^n$ and the unit ball $B^n$. These results will be utilized in our subsequent proofs.
	\begin{theorem}\cite[Theorem 1]{slot2023bpop} \label{binary_rate}
		Consider the problem \eqref{pop} with $d \leq n$ and $\mathit{X}= \bB^n$. Let $\xi_r^n$ be the least root of the degree-$r$ Krawtchouk polynomial with parameter $n$. Then if $(r+1)/n \leq 1/2$ and $d(d+1)\cdot \xi_{r+1}^n\leq 1/2$, we have:
		\begin{displaymath}
		    \dfrac{f_{\min}^X-\lb(f,\mathcal{Q}(\mathbb{B}^n))_r}{\|f\|_{\mathbb{B}^n}} \leq 2C_d\cdot \xi_{r+1}^n/n.
		\end{displaymath}
		Here $C_d > 0$ is a constant depending only on $d$, and  
  $\|f\|_{\mathbb{B}^n} := \max_{x \in \mathbb{B}^n}|f(x)|$. 
	\end{theorem}
	\begin{theorem}\cite[Theorem 3]{slot2111sum}\label{unit ball}
		Consider the problem \eqref{pop} with $\mathit{X}= B^n$. Then for any $r \geq 2nd$, the lower bound $\lb(f,\mathcal{Q}(B^n))_r$ for the minimal value $f_{\min}^X$  satisfies:
		\begin{displaymath}
		0\leq    f_{\min}^X-\lb(f,\mathcal{Q}(B^n))_r \leq \frac{C(n,d)}{r^2}
		\cdot(f_{\max}^X-f_{\min}^X).
		\end{displaymath}
		Here, $f_{\max}^X$ is the maximum value of $f$ over $B^n$, and $C(n,d)$ is a constant depending only on $n$ and $d$. In addition, this constant depends polynomially on $n$ (for fixed $d$) and polynomially on $d$ (for fixed $n$).
	\end{theorem}

    For a general feasible set, given a positive polynomial $f$ over such a set, numerous works have attempted to bound the order $k$ of the quadratic module such that $f$ belongs to 
    $\mathcal{Q}(X)_k$. Nie and Schweighofer \cite{nie2007complexity} established exponential bounds for this order in \cite[Theorem 8]{nie2007complexity}. More recently, Baldi and Mourrain \cite{baldi2021moment} proved polynomial bounds 
    based on the  Łojasiewicz inequality which we state below for the 
    set $\cX$ in \eqref{semialg by semid} but is also applicable to the set $X$ in 
    \eqref{semialg by poly}.
   
 \begin{theorem}\cite[Theorem 4.1]{dinh2016lojasiewicz}
 \label{Lojasiewicz}
For any compact set $K$ containing $\cX$ as defined  in \eqref{semialg by semid}, there exist a Łojasiewicz constant $C>0$ and a Łojasiewicz exponent $0 < L \leq 1$ depending on $\mathcal{X}$ such that
\begin{displaymath}
    d(x,\mathcal{X})\leq C\max\{0,-\lambda_m(x)\}^{L}\quad \forall x \in K.
\end{displaymath}
When $G(x)$ is a diagonal matrix 
consisting of $g_1(x),\ldots,g_m(x)$ on its diagonal,  
$\cX$ reduces to the set $X$ in 
 \eqref{semialg by poly}, and we get 
 \begin{displaymath}
     d(x,X)\leq C\max\big\{0,-\min_{j\in[m]}\{g_j(x)\} \big\}^{L}\quad \forall x \in K.
 \end{displaymath}
\end{theorem}
\begin{remark}
Note that the condition $L \leq 1$ is not present in \cite{dinh2016lojasiewicz}. However, we can assume it to be at most $1$ because of the compactness of $K$. In particular, if $L >1$, we can replace $L$ by $1$ and multiply the Łojasiewicz constant by $\max_{x \in K}\max\{0,-\lambda_m(x)\}^{L-1} <\infty$ to obtain a new inequality with the Łojasiewicz exponent $1$.
\end{remark}

With Theorem~\ref{Lojasiewicz}, we can now state 
the polynomial bounds provided by Baldi and Mourrain \cite{baldi2021moment}.
    \begin{theorem}\cite[Theorem 1.7]{baldi2021moment}\label{thm-Baldi}
        Let $\mathit{X}$ be defined as in \eqref{semialg by poly}, which is assumed to be contained in $[-1,1]^n$ and satisfied some normalization assumptions. Consider a positive polynomial $f$ over $\mathit{X}$ of degree $d$. Let $C,L$ be the Łojasiewicz coefficient and exponent given by Theorem 
   ~\ref{Lojasiewicz} with $\mathit{X}$ and $K=[-1,1]^n$. Then $f \in \mathcal{Q}(X)_k$ if
        \begin{align*}
 	      k \;\geq\; \gamma(n,X)d^{3.5n/L} E_f^{-2.5n/L},
 	  \end{align*}
 	  	where $\gamma(n,X) \geq 1$ depends only on $n$ and $X$. In the above, $E_f=f_{\min}^X/\|f\|_{[-1,1]^n}$ is a measure of how close $f$ is to having a zero on $\mathit{X}$.
    \end{theorem}
    
    The bounds wereimproved in the recent work \cite{baldi2025lojasiewicz}.
    \begin{theorem}\cite[Theorem 3.1]{baldi2025lojasiewicz}\label{thm-Baldi2}
        Let $\mathit{X}$ be defined as in \eqref{semialg by poly}, which is assumed to be contained in $B^n$ and satisfied some normalization assumptions. Consider a positive polynomial $f$ over $\mathit{X}$ of degree $d$. Let $C,L$ be the Łojasiewicz coefficient and exponent given by Theorem 
   ~\ref{Lojasiewicz} with $\mathit{X}$ and $K=B^n$. Then $f \in \mathcal{Q}(X)_k$ if
    \begin{displaymath}
        k = \mathcal{O}\left( n^2m\max_{j\in [m]}\{\deg g_j\}^6C^7E_f^{-7/L-3} \right),
    \end{displaymath}
    where $E_f=f_{\min}^X/\|f\|_{B^n}$.
    \end{theorem}
    For the matrix case, Huang \cite{huang2025complexity} extended Theorem~\ref{thm-Baldi2} to a set defined by a polynomial matrix inequality, as presented in \cite[Theorem 1.1]{huang2025complexity}. The corresponding bound is 
    $$
    k=\mathcal{O}\big(8^{7/L} 3^{6(m-1)}\theta(m)^3 n^2 {\rm deg}(G)^6
    C^7 d^{14/L} (f_{\min}/\|f\|_{B^n})^{-7/L-3}\big),
    $$
    where $\theta(m)\geq \frac{m! (m+1)!}{2^m}$.
    As we can observe, it depends at least exponentially on the matrix dimension $m$ of $G$. In contrast, we shall see later in Theorem~\ref{thm: bound on degree} that our bound depends 
    polynomially on $m$ (for a fixed $n$).

%While his original result applies to positive definite polynomial matrices, we focus here on reviewing its scalar version concerning positive polynomial.

%%%%%%%%%%%%%%%%%%%%%%%%%%

 %%%%%%%%%%%%%%%%%%%%%%%%%%%%%%%%%
	\subsection{Methodology}\label{sec: method}
	We introduce a family of polynomial penalty functions such that when added to the objective function of \eqref{SPOP}, produces a new polynomial optimization problem with two key properties:
    \begin{enumerate}
        \item The new optimal value closely approximates $f_{\min}$ in \eqref{SPOP}.
        \item The new optimal value can be approximated by solving a problem of the 
        form in \eqref{pop} over a {\em simple set} (namely, the binary 
        hypercube $\bB^n$, the unit ball $B^n$ or the hypercube $[-1,1]^n$).
    \end{enumerate}
    Leveraging Theorems~\ref{binary_rate} and~\ref{unit ball}, we analyze the convergence of this approach. Specifically, we consider \eqref{SPOP} in the form:
	\begin{align}\label{Spop}
		f_{\min} &= \min{f(x)}\nonumber\\
		& \mbox{subject to }  x\in \mathcal{X} = \{x \in\mathcal{C},\; G(x) \succeq 0\}.
	\end{align} 
	where $\cC$ is a simple set with known convergence rates for the standard scalar
	moment-SOS hierarchy. In this work, we consider two cases: $\mathcal{C}$ is the 
	binary hypercube in the discrete case (i.e., $\cX$ is a subset of $\bB^n$), and $\mathcal{C}$ is the $n$-dimensional unit ball in the continuous case. For convenience, we make the following assumption.
    \begin{assume}[\textbf{Normalization assumption}]\label{normalization assumption}
        The spectral radius of $G(x)$, denoted by $\rho(G(x))$, over $\cC$ is bounded by $1$, i.e., $\max_{x \in \cC}\rho(G(x)) \leq 1$.
    \end{assume}
 The above assumption can be satisfied in general by a proper scaling of $G$. 
    Since $\cC$ is a compact set, it is contained in some ball $B(0,R)$. Thus, the spectral radius of $G$ can be bounded as:
    \begin{displaymath}
       \rho(G(x)) \leq \sum_{1 \leq i,j \leq m}|g_{ij}(x)| \leq \sum_{1 \leq i,j \leq m}\|g_{ij}\|_{\infty}\max\{1,R^l\} \quad \forall \; x\in\cC.
   \end{displaymath}
    Therefore by scaling $G(x)$ by $\sum_{1 \leq i,j \leq m}\|g_{ij}\|_{\infty}\max\{1,R^l\}$, the normalization assumption is satisfied. Throughout this paper, we always specify the set $\cC$ before using this assumption. 
    
    We make a simple observation: If the optimal value $f_{\min}$ over $\mathcal{X}$ coincides with the optimal value of the objective function $f$ over $\mathcal{C}$, then the problem \eqref{Spop} is equivalent to the optimization problem over the simple set $\mathcal{C}$: 
	\begin{align*}
		f_{\min}^\cC &= \min\big\{{f(x)}:\
		x \in \mathcal{C}\big\},
	\end{align*} 
    where the convergence rate has been shown to be $\mathcal{O}(1/r^2)$ in either Theorem~\ref{binary_rate} or~\ref{unit ball}.
	When the optimal value $f_{\min}$ is larger than $f^\cC_{\min}$ (the optimal value of $f$ over $\mathcal{C}$), we employ a polynomial penalty function $P(x)$ satisfying that:
    \begin{eqnarray*}
        P(x) \approx 0 \; \forall\; x \in \cX, \quad 
        P(x) \gg \; 0\; \forall\; x \in \cC \backslash \cX.
    \end{eqnarray*}
    Then, we consider the penalized POP below:
	\begin{align} \label{penaltySpop}
		F_{\min} &= \min \big\{ F(x)\;:=\;f(x)+P(x) :\; x \in \mathcal{C}\big\}.
	\end{align}
    The key properties of $P(x)$ implies that $F_{\min} \approx f_{\min}$, and $F_{\min}$ can be approximated by the moment-SOS hierarchy over $\mathcal{C}$ with the convergence rate of $\mathcal{O}(1/r^2)$.
    Our proposed penalty function takes the form:
    \begin{displaymath}
        P(x)\;=\; -\langle h(G(x)),G(x) \rangle,
    \end{displaymath}
    where $h\in \cP_+([-1,1])$, which approximates a step function on $[-1,1]$ that is equal to $0$ in the interval $[0,1]$, and equal to a large positive number $N$ in $[-1,0)$. 
    By expressing $h$ as $h(t) = \sum_{i=0}^k h_i t^i$, and using the property of trace, we have that
    \begin{align}\label{eigenvalue equality}
        P(x)&\;= -\langle h(G(x)),G(x) \rangle= -\sum_{i=0}^kh_i\langle G(x)^i,G(x) \rangle
        \;=  -\sum_{i=0}^kh_i\langle G(x)^{i+1},I_m \rangle
    \nonumber \\
        &\;= -\sum_{i=0}^kh_i\tr(G(x)^{i+1})
        \;= -\sum_{i=0}^kh_i\sum_{j=1}^m\lambda_j(x)^{i+1}
        \;=-\sum_{j=1}^m\Big(\sum_{i=0}^kh_i\lambda_j(x)^{i+1}\Big)
        \nonumber \\
        & \;= -\sum_{j=1}^m\lambda_j(x)h(\lambda_j(x)).
    \end{align}
    Since $h(t)$ approximates the step function as we describe above, we observe the following properties:
    \begin{itemize}
        \item If $x \in \cX$, then $G(x) \succeq 0$ and all its eigenvalues $\lambda_j(x) \in [0,1]$ under  Assumption~\ref{normalization assumption}. Whence, $h(\lambda_j(x)) \approx 0 \;\forall j \in [m]$ 
        and \eqref{eigenvalue equality} implies that 
        $P(x) \approx 0$ for $x \in \cX$.
        \item If $x \in \cC \backslash \cX$, then $G(x)$ has all eigenvalues in $[-1,1]$ and at least one eigenvalue in $[-1,0)$. For those eigenvalue $\lambda_j(x) \geq 0$, we have that $\lambda_j(x)h(\lambda_j(x)) \approx 0$, and for those eigenvalue $\lambda_j(x) < 0$, $-\lambda_j(x)h(\lambda_j(x)) \gg 0$. Hence, $P(x) \gg 0$. 
    \end{itemize}
    Therefore, $P(x)$ intuitively satisfies the key properties of the penalty function we need. Moreover, the non-negativity of $h$ over $[-1,1]$ can be characterized by the Fekete, Markov-Luk\`acz theorem \cite[Theorem 3.1]{deklerk2017improved}, which leads to a nice construction of our novel hierarchy without increasing the size of the SDP matrix variable by 
    a factor of $m$ as described in Section~\ref{sec: summary on complexity}. The mathematical arguments and details are elaborated in the the next two sections.

%%%%%%%%%%%%%

    %%%%%%%%%%%%%%%%%%%%%%%%%%%%%%%%%%%%%
\section{Binary polynomial optimization problems with polynomial matrix semidefinite constraints}\label{section: BPOP}

 \subsection{Matrix moment-SOS hierarchy}
     Consider \eqref{SPOP} in the case $\cX$ where is a subset of the binary hypercube $\bB^n$, yielding the following formulation:
	\begin{align}\label{BPOP}
		f_{\min} &= \min{f(x)}\nonumber \\
		& \mbox{subject to } x \in \mathcal{X}:=\{x \in \mathbb{B}^n:\ G(x) \succeq 0\}.
	\end{align}
  Notice that any polynomial objective function over any subset of $\bB^n$ can be replaced by its image in the quotient ring $\mathbb{R}[x]/\langle x_1-x_1^2,\dots,x_n-x_n^2 \rangle$, where all polynomials have degree at most $n$. Thus, we assume that the degree $d$ of $f$ is at most $n$ without loss of generality. We propose the following hierarchy of lower bounds for 
  \eqref{BPOP}:
    \begin{equation*}
            \lb(f, \widetilde{\cQ}(\cX))_r :=\; \max\{t \in \bR:\; f-t \in  \widetilde{\cQ}(\cX)_{2r}\},
    \end{equation*}
    where $\widetilde{\cQ}(\cX)_{2r}=\; \cQ(\bB^n)_{2r}+\cH(\cX)_{\lfloor r/l\rfloor-1}$ is defined as in \eqref{bpop-qmodule}.
    
    \begin{proposition}\label{prop: SDP certificate}
        For any non-negative integer $k$, verifying whether  
       a polynomial $p(x) \in \bR[x]_{2l(k+1)}$ is contained in 
        $\cH(\cX)_{k}$ can be verified by an SDP with
         $(\lceil k/2 \rceil+1)^2$ variables.
    \end{proposition}
    \begin{proof}
        For any non-negative univariate polynomial $h$ over $[-1,1]$ of degree at most $2k$, by the Fekete, Markov-Luk\`acz theorem \cite[Theorem 3.1]{deklerk2017improved}, there exist SOS univariate polynomials $h_1$ and $h_2$ of degree at most $2\lceil k/2 \rceil$ and $2\lceil k/2 \rceil-2$, respectively, that satisfy the following representation:
        \begin{displaymath}
            h(t) \;=\; h_1(t)+(1-t^2)h_2(t).
        \end{displaymath}
    Therefore, $\langle h(G(x)),G(x) \rangle$ can be expressed as:
    \begin{displaymath}
        \langle h(G(x)),G(x) \rangle \;=\; \langle h_1(G(x)),G(x) \rangle+ \langle h_2(G(x)),(I_m-G(x)^2)G(x) \rangle.
    \end{displaymath}
    By writing $h_1$ and $h_2$ explicitly in terms of the standard monomials: 
    \begin{displaymath}
         h_1(t)=\; \sum_{i=0}^{2\lceil k/2 \rceil}h_i^{(1)}t^i,\quad h_2(t)= \sum_{i=0}^{2\lceil k/2 \rceil-2}h_i^{(2)}t^i,
    \end{displaymath}
    we can express $\langle h(G(x)),G(x) \rangle$ as in \eqref{eigenvalue equality} 
    using the trace of powers of $G(x)$:
    \begin{equation*}\label{prop: SDP size eq 1}
        \langle h(G(x)),G(x) \rangle =\; \sum_{i=0}^{2\lceil k/2 \rceil}h_i^{(1)}\tr(G(x)^{i+1})+\sum_{i=0}^{2\lceil k/2 \rceil-2}h_i^{(2)}(\tr(G(x)^{i+1})-\tr(G(x)^{i+3}).
    \end{equation*}
    Let $b_{\lceil k/2 \rceil}(t)$ be the standard monomials basis $(1,\; t,\;\ldots,\;t^{\lceil k/2 \rceil})^{\top}$. As in \cite{lasserre2009moments}, we say that $h_1$ and $h_2$ are SOS polynomials if there exist positive semidefinite matrices $H_1$ and $H_2$
     with rows and columns indexed  respectively by the basis $b_{\lceil k/2 \rceil}(t)$ and $b_{\lceil k/2 \rceil-1}(t)$ satisfying that 
    \begin{equation}\label{prop: SDP size eq 2}
        h_1(t)=\; \langle H_1, b_{\lceil k/2 \rceil}(t)b_{\lceil k/2 \rceil}(t)^{\top} \rangle,\quad h_2(t)=\;
        \langle H_2, b_{\lceil k/2 \rceil-1}(t)b_{\lceil k/2 \rceil-1}(t)^{\top} \rangle.
    \end{equation}
    For any $v \in \bN$, we define the following polynomial matrices with rows and columns indexed by the basis $b_v(t)$ as follows: for $0 \leq i,j \leq v$, their $(i,j)$-element are given as
    \begin{equation}\label{def: trace matrix}
        [P_v^G(x)]_{i,j}=\; \tr(G(x)^{i+j+1}),\quad  [Q_v^G(x)]_{i,j}=\; \tr(G(x)^{i+j+1})-\tr(G(x)^{i+j+3}).
    \end{equation}
    We next plug \eqref{def: trace matrix} into \eqref{prop: SDP size eq 2} to obtain the following representation:
    \begin{equation}\label{prop: SDP size eq 3}
        \langle h(G(x)),G(x) \rangle =\; \langle H_1,P_{\lceil k/2 \rceil}^G(x)\rangle+\langle H_2, Q_{\lceil k/2 \rceil-1}^G(x) \rangle.
    \end{equation}
    We note that the degree of both $\langle H_1,P_{\lceil k/2 \rceil}^G(x)\rangle$ and $\langle H_2, Q_{\lceil k/2 \rceil-1}^G(x) \rangle$ are at most $\deg(G)(2\lceil k/2 \rceil+1)\leq 2l(k+1)$.
    As a result, for any $p(x) \in \bR[x]_{2l(k+1)}$, the problem of checking the membership of $p$ in $\cH(\cX)_k$ can be cast as the following semidefinite feasibility problem: find
    $H_1$ and $H_2$ such that  
    \begin{equation}\label{membership of H}
        p(x) =\; \langle H_1,P_{\lceil k/2 \rceil}^G(x)\rangle+\langle H_2, Q_{\lceil k/2 \rceil-1}^G(x) \rangle,\quad H_1 \in \bS_+^{\lceil k/2 \rceil+1},\quad H_2 \in \bS_+^{\lceil k/2 \rceil}.
    \end{equation}
    The number of variables in this SDP is calculated as the sum of the variables in the symmetric matrices $H_1$ and $H_2$ as follows:
    \begin{displaymath}
        (\lceil k/2 \rceil+1)(\lceil k/2 \rceil+2)/2
        +\lceil k/2 \rceil(\lceil k/2 \rceil+1)/2 = (\lceil k/2 \rceil+1)^2.
    \end{displaymath}
    \end{proof}
    \begin{corollary}\label{cor: proposed momemnt SDP for bpop}
    For any integer $r\geq 3l$,
        the dual SDP problem corresponding to the hierarchy \eqref{proposed hierarchy for bpop} is 
        \begin{align}\label{proposed momemnt SDP for bpop}
            \min_{y\in \bR^{s(n,2r)}}&\quad L_y(f)\\
            \mbox{subject to}&\quad y_0 =1,\; \M_r(y) \succeq 0,\; \M_{r-1}((x_i(x_i-1))y)=0\; \forall i \in [n],\nonumber\\
            &\quad L_y(P_{\lceil(\lfloor r/l \rfloor -1)/2 \rceil}^G)\succeq 0,\; 
            L_y(Q_{\lceil(\lfloor r/l \rfloor -1)/2 \rceil-1}^G) \succeq 0.\nonumber
        \end{align}
    \end{corollary}
    \begin{proof}
        We can express $b_r(x)b_r(x)^{\top} = \sum_{\alpha \in \bN^n_{2r}}B_{\alpha}x^{\alpha}$ for suitable coefficient matrices $B_\alpha$ with $\alpha \in \bN^n_{2r}.$
        For $i \in [n]$, we can write $x_i(x_i-1)b_{r-1}(x)b_{r-1}(x)^{\top} = \sum_{\alpha \in \bN^n_{2r}}B^{(i)}_{\alpha}x^{\alpha}$. Let
        $k=\lfloor r/l \rfloor -1$, we can similarly write 
        $P_{\lceil k/2 \rceil}^G(x)= \sum_{\alpha \in \bN^n_{2r}}P_{\alpha}x^{\alpha}$, 
        $Q_{\lceil k/2 \rceil-1}^G(x)= \sum_{\alpha \in \bN^n_{2r}}Q_{\alpha}x^{\alpha}$ for suitable coefficient matrices $P_\alpha$ and $Q_\alpha$
        with $\alpha \in \bN^n_{2r}.$ Then \eqref{proposed hierarchy for bpop} can be cast as the following SDP: 
        \begin{align}\label{SOS SDP bpop}
            \max &\ t\\
            \mbox{subject to:}&\ f_{\alpha}= \langle X, B_{\alpha} \rangle+\sum_{i=1}^n\langle X_i, B^{(i)}_{\alpha} \rangle+\langle Y, P_{\alpha} \rangle + \langle Z, Q_{\alpha} \rangle+t \delta_{\{\alpha=0\}}, \;\forall\;\alpha \in  \bN^n_{2r}, 
            \nonumber\\
            &\ X \in \bS_+^{s(n,r)},\; X_i \in \bS_+^{s(n,r-1)}\; \forall i \in [n],\quad
             Y \in \bS_+^{\lceil k/2 \rceil+1},\; Z \in \bS_+^{\lceil k/2 \rceil}.\nonumber
        \end{align}
        Then by the duality theory of SDP, the dual of the SDP problem corresponding to \eqref{proposed hierarchy for bpop} is given by \eqref{proposed momemnt SDP for bpop}.
    \end{proof}
    \begin{remark}
        Under the Archimedean condition, the strong duality holds between \eqref{proposed hierarchy for bpop} and \eqref{proposed momemnt SDP for bpop} (see e.g., \cite{josz2016strong}). In comparison with the existing SDP relaxation \eqref{primal SDP for SPOP}, we observe  that the condition $\M_{r-l}(Gy) \succeq 0$ has been replaced by the conditions
        \begin{displaymath}
            L_y(P_{\lceil(\lfloor r/l \rfloor -1)/2 \rceil}^G)\succeq 0,\quad  L_y(Q_{\lceil(\lfloor r/l \rfloor -1)/2 \rceil-1}^G) \succeq 0,
        \end{displaymath}
        whose matrix sizes are at most $\lceil(\lfloor r/l \rfloor -1)/2 \rceil$, which is independent of $m$ and $n$, and is
        significantly smaller than $m \cdot s(n,r-l)$, which is the size of the constraint matrix $\M_{r-l}(Gy)$ in Hol-Scheiderer's Hierarchy \eqref{primal SDP for SPOP}. As a result, we have successfully avoided the explosion of the matrix size in the SDP relaxation. In what follows, we study the 
        convergence rate of \eqref{proposed momemnt SDP for bpop} as $r$ increases.
    \end{remark}
 %%%%%%%%%%%%%%%%%%%%%%%%%%%%%%
	\subsection{Simple representations of $\mathcal{X}$}
    In this section, we show that the polynomial matrix inequality in the description of $\cX$:
    \begin{displaymath}
        \cX =\; \{ x \in \bB^n\;:\; G(x) \succeq 0\}
    \end{displaymath}
    can be equivalently expressed through a scalar polynomial inequality of the form $\langle h(G(x)),G(x) \rangle \geq 0$. In the next theorem, we demonstrate this replacement property to arbitrary discrete sets beyond the binary hypercube $\bB^n$. This not only shows the key insight of using members of $\cH(\cX)$ to effectively reduce the matrix inequality $G(x) \succeq 0$ to a scalar condition, but also provides a glimpse of the 
    construction of the penalty function framework in Section~\ref{sec: method}.
	  \begin{theorem}\label{simple_rep2} Consider the following set $\cX$ defined by
      \begin{displaymath}
        \cX =\; \{ x \in \cC\;:\; G(x) \succeq 0\}
    \end{displaymath}
    that satisfies Assumption \ref{normalization assumption} with $\cC=\bB^n$. Then there exists a non-negative univariate polynomial $h(t) \in \mathbb{R}[t]$ over the interval $[-1,1]$ such that
	  	\begin{displaymath}
	  	    \mathcal{X} = \{ x \in \cC: \langle h(G(x)),G(x) \rangle \geq 0 \}. 
	  	\end{displaymath} 
	  \end{theorem}
      
	\begin{proof}
	For any $ x \in \mathcal{X}$, Assumption \ref{normalization assumption} implies that $\lambda_i(x) \in [0,1] \; \forall i \in [m]$. Therefore, for any non-negative univariate polynomial $h(t)$ over $[-1,1]$, the following inequalities holds
    \begin{displaymath}
        \langle h(G(x)),G(x) \rangle =\; \sum_{i=1}^n \lambda_i(x)h(\lambda_i(x)) \geq\; 0.
    \end{displaymath}
    This inequality induces the following containment:
	\begin{displaymath}
	  \mathcal{X} \subset \{ x \in \mathcal{C}:\; \langle h(G(x)),G(x) \rangle \geq 0 \}.
	\end{displaymath}
    
	We next show the reverse containment for some $h \in \cP_+([-1,1])$. Due to the discreteness of $\cC=\bB^n$ and the normalization of $G(x)$ in 
    Assumption \ref{normalization assumption}, there exists a negative number $\lambda$ such that $0 > \lambda := \max_{x \in \cC\backslash \mathcal{X}}\lambda_m(x) \geq -1$. 
    
    We observe that for any $t \in [-1, \lambda]$, the function $h(t)$ needs to be much larger than the value of $h(t)$ over the interval $[0,1]$ so that if $G(x)$ has some negative eigenvalue $\lambda_j(x)$, the function value $\lambda_j(x)h(\lambda_j(x))$ will dominate that of the other non-negative eigenvalues to make $\langle h(G(x)),G(x) \rangle$ negative. This behavior can be achieved by a non-negative polynomial approximation of the piecewise affine function $a(t)$ defined as:
	  	\begin{displaymath}
	  	    a(t) = \begin{cases}
	  		   1 &\mbox{if} \ t \in [0,1]\\
	  		   -m/\lambda  & \mbox{if}\ t \in [-1,\lambda]\\ 
	  		   1-\Big(1+\dfrac{m}{\lambda}\Big)\cdot \dfrac{t}{\lambda}& \mbox{if}\ t \in [\lambda,0].
	  	\end{cases} 
        \vspace{-1em}
	  	\end{displaymath}\begin{figure}[tbhp]
        \centering
        \begin{subfigure}[t]{0.45\textwidth}
            \centering
            \includegraphics[width=\textwidth]{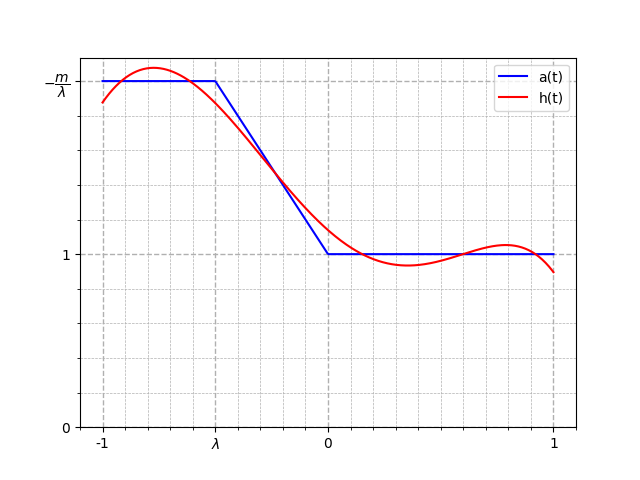}
            
        \end{subfigure}
        \hfill
        \begin{subfigure}[t]{0.45\textwidth}
            \centering
            \includegraphics[width=\textwidth]{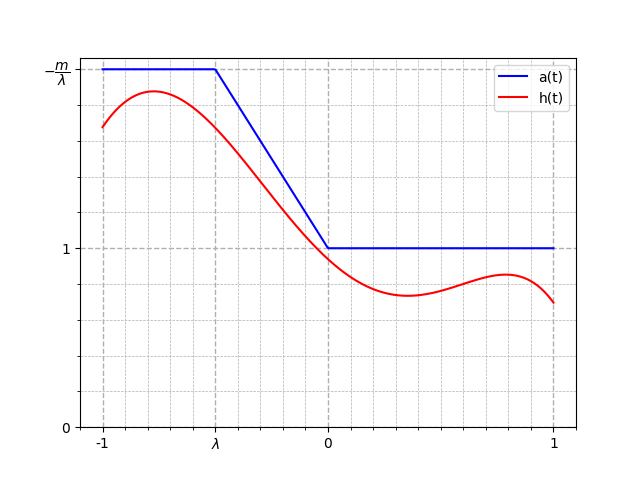}
            
        \end{subfigure}
        \caption{The left panel shows the 
        polynomial approximation of $a(t)$, and the right panel 
        shows the polynomial approximation of $a(t)$ from below after a 
        vertical translation.}
        \label{fig:approximation}
        \end{figure}

         Given a positive $\varepsilon >0$, the Stone--Weierstrass theorem (see \cite{rudin1953principles}) ensures that there exists a polynomial $h(t)$ such that $\| h- a\|_{[-1,1]} \leq \varepsilon/2$. We then subtract $\varepsilon/2$ from $h(t)$ to obtain a polynomial approximation of $a(t)$ from below, which is illustrated in Figure~\ref{fig:approximation}. Therefore, we can choose $h(t)$ to be a polynomial approximation of $a(t)$ from below  satisfying that $\| h- a\|_{[-1,1]} \leq \varepsilon$.
        
       Let us choose $\varepsilon <1/2$. Since $a(t) \geq  1 \ \forall t \in [-1,1]$, $ h(t) \geq a(t)-\varepsilon \geq 1- \varepsilon > 0 \ \forall t \in [-1,1]$. Moreover, for any $t \in [-1,\lambda]$, $h(t) \geq a(t) -\frac{1}{2}= \frac{-m}{\lambda}-\frac{1}{2} \geq m-\frac{1}{2} >0$. Thus, the following properties hold:
\begin{eqnarray*}
&\forall \; t \in [0,1], \quad 0<h(t) \leq a(t) =1 &\quad \Rightarrow \quad 0 \;\leq\;  th(t) \;\leq\; 1,\\
&\forall\;  t \in [-1,\lambda], \quad h(t) \geq \dfrac{-m}{\lambda}- \dfrac{1}{2} > 0 &\quad \Rightarrow \quad th(t) \;\leq\;\left(\dfrac{-m}{\lambda}-\dfrac{1}{2} \right)  \lambda \;<\; 0.
\end{eqnarray*}

       For any $x \in \mathcal{C}\backslash\mathcal{X}$, $G(x)$ is not a positive semidefinite matrix. Then by Assumption \ref{normalization assumption}, we obtain the followings:
       \begin{itemize}
           \item The smallest eigenvalue $\lambda_m(x) \in [-1,\lambda]$, and $h(\lambda_m(x))\lambda_m(x) \leq \left(\frac{-m}{\lambda}-\frac{1}{2} \right)\lambda$.
           \item For any $i \in [m-1]$:
           \begin{enumerate}
               \item[(i)] If $\lambda_i(x)$ is non-negative, then $\lambda_i(x) \in [0,1]$, and $\lambda_i(x)h(\lambda_i(x)) \leq 1$;
               \item[(ii)] If $\lambda_i(x)$ is negative, then $\lambda_i(x)h(\lambda_i(x)) <0 < 1$ since $h(t)> 0 \; \forall t \in [-1,1]$.
           \end{enumerate}
       \end{itemize}
       This yields the key inequality:
	  \begin{align*}
	      \langle h(G(x)),G(x) \rangle = h(\lambda_m(x))\lambda_m(x) + \sum_{i=1}^{m-1}h(\lambda_i(x))\lambda_i(x)
            \leq  
            \Big(\dfrac{-m}{\lambda}-\dfrac{1}{2}\Big)\lambda +m-1 <0.
	  \end{align*}
      Thus $x\in\cC\backslash\cX \Rightarrow \langle h(G(x)),G(x) \rangle <0$. 
Hence $ \left\{x \in \mathcal{C} : \langle h(G(x)),G(x) \rangle \geq 0 \right\}\subset\cX$.
	  \end{proof}

      \begin{comment}
      \begin{corollary}\label{thm: convergence for bpop}
          Under Assumption~\ref{normalization assumption} with $\cC =\bB^n$, if $f$ is a positive polynomial over $\cX$ defined as in \eqref{BPOP}, then $f \in \overline{\cQ}(\cX)$.
      \end{corollary}
      \begin{proof}
          We recall the construction of the univariate polynomial $h$ in the proof of Theorem~\ref{simple_rep2}, it is still valid if we replace $1$ and $m/\lambda$ by any positive numbers. Whence, there exists a positive univariate polynomial $h$ over $[-1,1]$ such that 
          \begin{displaymath}
              \langle h(G(x)),G(x) \rangle < \min_{x \in \cX} f \ \forall x \in \cX, \text{ and } \langle h(G(x)),G(x) \rangle < -\|f\|_{\bB^n} \ \forall x \in \bB^n \backslash \cX.
          \end{displaymath}
          Therefore, $f(x) - \langle h(G(x)),G(x) \rangle $ is positive over $\bB^n$, then 
          \begin{displaymath}
              f(x) - \langle h(G(x)),G(x) \rangle \in \cQ(\bB^n) \; \Rightarrow \; f \in \overline{\cQ}(\cX).
          \end{displaymath}
      \end{proof}
      \end{comment}
      
   %%%%%%%%%%%%%%%%%%%%%%%%%%%%%%%%%%
	  \subsection{Approximation of a nonnegative piecewise affine function by 
	  nonnegative polynomials} \label{section: Penalty for binary}
	  
	  Let $N$ be a given positive number and $\lambda \in [-1,0)$ be a given negative number.
	  Consider the following nonnegative piecewise affine function defined over the interval $[-1,1]$:
	  \begin{eqnarray}
	   a(t) = \begin{cases}
	  		0 &\mbox{if} \ t \in [0,1],\\
	  		N  & \mbox{if}\ t \in [-1,\lambda],\\ 
	  		N t/\lambda  & \mbox{if}\ t \in [\lambda,0].
	  	\end{cases}
	 \end{eqnarray} 
	  Our goal is to construct a non-negative univariate polynomial $h(t)$ to approximate 
	  $a(t)$ over $t\in [-1,1]$. According to Jackson's theorem, the approximation error is controlled by the smoothness 
	  of the underlying function to be approximated. Thus instead of directly
	  approximating $a(t)$ by a polynomial, for a given positive integer $k$, we first 
	  approximate $a(t)$ by a $k$ times differentiable nonnegative function $q(t)$ defined by:
      \begin{equation}\label{eq-q1}
	  		    q(t) = \begin{cases}
	  		0 &\mbox{if} \ t \in [0,1],\\
	  		N  & \mbox{if}\ t \in [-1,\lambda],\\ 
	  		c(t)& \mbox{if}\ t \in [\lambda,0],
	  	\end{cases}
	  		\end{equation}
         where  
         $c(t)$ is a {\em concatenation polynomial} chosen so that $q(t)$ is $k$ times continuously differentiable on $[-1,1]$. We then approximate $q(t)$ from below using a Chebyshev polynomial $h$ of degree $v$, with the error controlled by Jackson's theorem. 
         In order for the piecewise polynomial function $q(t)$ to be $k$ times differentiable, 
        we look for a concatenation polynomial $c(t)$ that satisfies the following conditions:
	  	\begin{equation}\label{eq-c}
	  	\left\{ \begin{array}{l}
	  		\mbox{$c(0) = 0$, \ $c(\lambda) = N$, \ $c(t) \geq 0$  for all $t\in [\lambda,0]$,} \\
	  		\mbox{$c^{(i)}(0)= c^{(i)}(\lambda) =0$, for all $i \in [k]$,}
	  		\end{array}
	  		\right.
	  	\end{equation}
	  	where $c^{(i)}(\cdot)$ denotes the $i$-th derivative of $c(\cdot)$. Then, 
	  	since the function $q(t)$ is a piece-wise polynomial, the one-sided derivatives of any order exist, and the above conditions ensure that $q$ is continuously differentiable up to the order $k$.

	  We can first construct $c(t)$ on the interval $[0,1]$ and then scaling $t$ by $\lambda$ to 
	  fit the interval $[\lambda,0]$. Thus, we can replace the conditions in \eqref{eq-c} by
	  \begin{equation}\label{eq-c0}
	  	\left\{ \begin{array}{l}
	  		\mbox{$c(0) = 0$,\ $c(1) = 1$,\ $c(t) \geq 0$  for all $t\in [0,1]$,} \\
	  		\mbox{$c^{(i)}(0)= c^{(i)}(1) =0$, for all $i \in [k]$.}
	  		\end{array}
	  		\right.
	  	\end{equation}
   	  A simple observation is that the polynomial $c(t)-t$ behaves similarly at $t=0$ and $t=1$. Thus, the symmetry of $c(t)-t$ can reduce the number of conditions on $c(t)$ as shown in the following lemma.
      
	  \begin{lemma}\label{concatination conditon}
	  	Consider $c(t)= t +(2t-1)T(t(1-t))$, where $T(t)$ is a univariate polynomial with constant coefficient $0$.  Then $c(0) = 0$, $c(1) =1$.
    Moreover, if $ c^{(i)}(0)=0 \ \forall i \in [k]$, then $ c^{(i)}(1)=0 \ \forall i\in [k]$.
	  \end{lemma}
	  \begin{proof}
	  	Since $c(t) = t +(2t-1)T(t(1-t))$, and the constant coefficient of the polynomial $T(t)$ is $0$, we obtain that  $c(0) = 0$ and $c(1) = 1$. Moreover, the following equality holds true for any $t \in [0,1]$:
	  	\begin{displaymath}
	  	    c(t) + c(1-t) = t + 1-t +(2t-1)T(t(1-t))+(1-2t)T((1-t)t)=1.
	  	\end{displaymath}
	  	Therefore, by the chain rule, for any $i \in [k]$, we obtain that
	  	\begin{displaymath}
	  	    c^{(i)}(t)+(-1)^ic^{(i)}(1-t)=0.
	  	\end{displaymath}
        Thus $c^{(i)}(0)=0$ also implies that $c^{(i)}(1)=0\ \forall i \in [k]$.
	  \end{proof}
%   \begin{remark}
%       We state here the reasons why we choose $c(t) = t +(2t-1)T(t(1-t))$. First, we observe that $c(0)=0$ and $c(1)=1$. Then, subtracting $t$ from $c(t)$ leads to the fact that $t = 0$ and $t=1$ are the roots of $c(t)-t$. Second, the symmetry of $c(t)$, which helps to reduce the number of conditions on $c(t)$, also implies that $t = 1/2$ is a root of $c(t)-t$. So $c(t) -t $ should be expressed in terms of a symmetric function around $1/2$ on $[0,1]$ vanishing at $t=0$ and $t=1$, such as $t(1-t)$. Hence, we propose the form $c(t) = t+(2t-1)T(t(1-t))$.
%   \end{remark}
   
   Next, we discuss how to construct the polynomial $T(t)$ in Lemma \ref{concatination conditon}.
	Since the graph of $c(t) = t +(2t-1)T(t(1-t))$ for $t\in[0,1]$ is symmetric about the point $(1/2,1/2)$, we only need to focus on the interval $[0,1/2]$. Set $u(t) = t(1-t)$ for $t\in (-\infty,1/2]$. Then $u$ is a one-to-one function mapping $(-\infty, 1/2]$ to $(-\infty,1/4]$, and the inverse function is $t(u) = \frac{1}{2}(1-\sqrt{1-4u})$ for all $u \in (-\infty,1/4]$. We define 
      \begin{displaymath}
          \varphi(u) =\; \frac{2u}{1-4u + \sqrt{1-4u}},\quad u \in (-\infty,1/4).
      \end{displaymath}
      It is clear that 
      \begin{eqnarray*}
          \varphi\circ u (t) %=\; \frac{2u}{1-4u + \sqrt{1-4u}}
          =\; \frac{2t(1-t)}{1-4t(1-t)+\sqrt{1-4t(1-t)}}
          %\\=\; \frac{2t(1-t)}{(1-2t)^2+(1-2t)}
          =\; \frac{t}{1-2t}, \quad \forall t \in (-\infty,1/2).
      \end{eqnarray*}
    Therefore, $t+ (2t-1)\varphi (u(t)) =0$ for $t \in (-\infty,1/2)$. We then choose $T$ to be the $k$-th order Taylor polynomial of $\varphi(u)$ at $u=0$, which is proved to be our desired polynomial in Lemma~\ref{lem: Taylor expansion of varphi}. 
    In order not to break the flow of the current presentation, we postpone 
    the proofs related to $\varphi(u)$ to Appendix~\ref{appendix: univariate polynomial}.    
    With the choice of $T$ in Lemma~\ref{lem: Taylor expansion of varphi}, we can construct the required concatenation polynomial
    satisfying the required conditions \eqref{eq-c0}.
  
	  \begin{proposition}\label{prop: concatenation poly}
	  	Let $T_k(u) = \sum_{i=0}^k a_i u^i$ be the Taylor polynomial of order $k\geq 1$ of $\varphi(u)$ at $u=0$. Then it induces a polynomial of degree at most $2k+1$ defined as 
	  	\begin{displaymath}
	  	    c_k(t)=t+(2t-1)T_k(t(1-t)),
	  	\end{displaymath}
	  	which satisfies the following properties
	  	\begin{eqnarray}
	  	\left\{
	  	\begin{array}{l}
	  		\mbox{$c_k(0)=0$,\ $c_k(1)=1$,\ $ c_k^{(i)}(0)=c_k^{(i)}(1)=0\ \forall i \in [k]$,}
	  		\\[3pt]
	  		\mbox{$0\leq c_k(t) \leq 1$,\ $c_k(t)+c_k(1-t) =1 \ \forall \;t \in [0,1]$.}
	  		\end{array} \right.
	  	\end{eqnarray} 
	  \end{proposition}
	 
\begin{figure}[tbhp]
        \centering
        \includegraphics[width=0.45\linewidth]{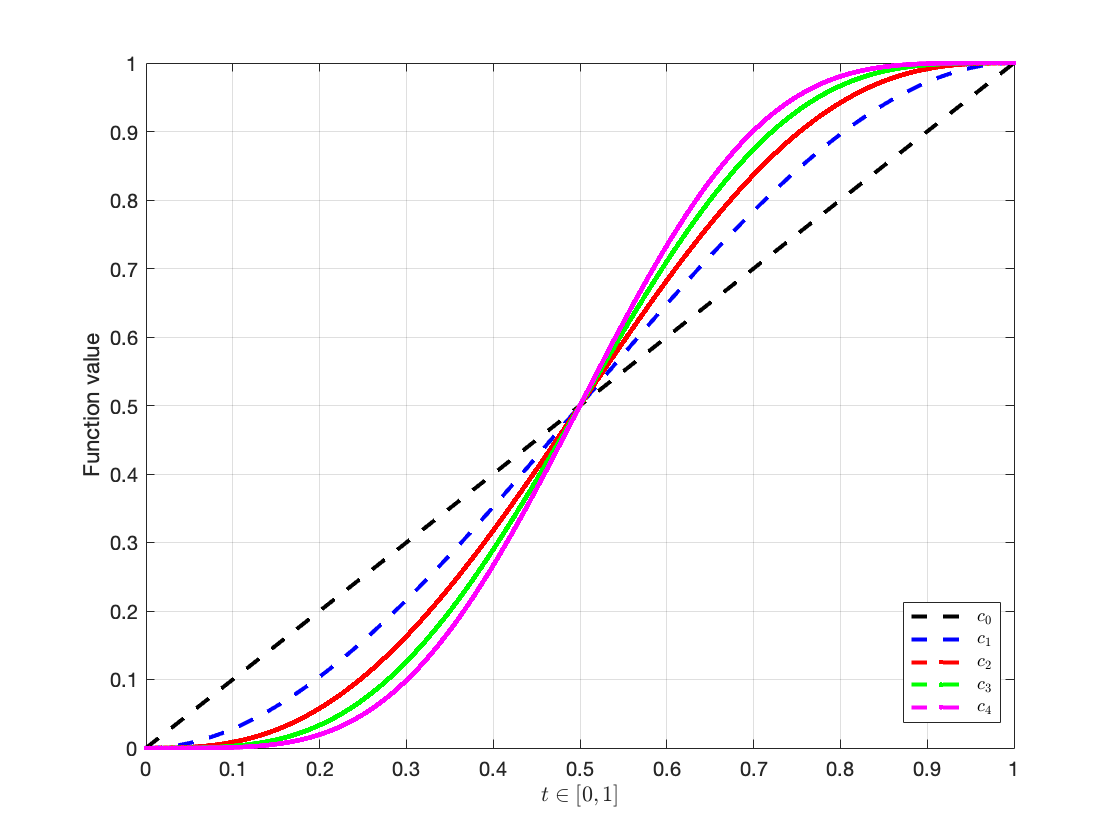}
        \caption{Comparison of $c_0, c_1,c_2,c_3 $ and $c_4$ over the interval $[0,1]$ in terms of their smoothness at the end points of the interval.}
        \label{fig:main}
        \end{figure}

	  \begin{remark} We show the first few concatenation polynomials as follows:
	  	\begin{eqnarray*}
	  	&&	c_0(t) = t,\quad c_1(t) =  t +(2t-1)u =  t +(2t-1)t(1-t), \\
	  	&&	c_2(t) = t + (2t-1)(u +3u^2) =  t +(2t-1)\big(t(1-t)+3t^2(1-t)^2\big),\\
	  	&&	c_3(t) = t +(2t-1)(u +3u^2+10u^3) = t +(2t-1)\big(t(1-t)+3t^2(1-t)^2+10t^3(1-t)^3\big).
%	  		\item $c_4(t) = t +(2t-1)t(1-t)+3(2t-1)t^2(1-t)^2+10(2t-1)t^3(1-t)^3+35(2t-1)t^4(t-1)^4=t+(2t-1)(u +3u^2+10u^3+35u^4)$.
	  	\end{eqnarray*}
	  	    Figure \eqref{fig:main} illustrates that $c_k(\cdot)$ is $k$ times differentiable at the end points of $[0,1]$.
	  \end{remark}
  \vspace{-1em}

        With the polynomial $c_k(t)$ in Proposition~\ref{prop: concatenation poly}, we can construct the following candidate function for
        $q(t)$ in \eqref{eq-q1}: 
	  \begin{equation}\label{smooth penalty}
	      q_k(\lambda,N)(t)=\begin{cases}
	  	0 &\mbox{if} \ t \in [0,1]\\
	  	N  & \mbox{if}\ t \in [-1,\lambda]\\ 
	  	Nc_k(t/\lambda)& \mbox{if}\ t \in [\lambda,0].
	  \end{cases}
	  \end{equation}
   According to Proposition~\ref{prop: concatenation poly}, $c_k^{(i)}(0) = c_k^{(i)}(1)=0 \ \forall i \in [k]$. Therefore, $q_k(\lambda,N)$ is $k$ times continuously differentiable 
   and non-negative over $[-1,1]$. Moreover, the $k$-th derivative of 
   $q_k(\lambda,N)$ has bounded variation on $[-1,1]$ as stated in the following lemma. 
   For simplicity, we recall the definitions of absolute continuity, total variation and functions of bounded variation, along with the related proofs, in Appendix~\ref{appendix: univariate polynomial}.
   
   \begin{lemma}\label{lemma: AC and BV}
       For any positive integer $k$, the $(k-1)$-th derivative of $q_k(\lambda,N)$ is absolutely continuous, and the $k$-th derivative of $q_k(\lambda,N)$ is of bounded variation on $[-1,1]$ with the total variation bounded by $V=3N|\lambda|^{-k}4^{k}k!k$.
   \end{lemma}
   \begin{proof}
       The proof is given in Appendix~\ref{appendix: univariate polynomial}, where all related definitions are stated.
   \end{proof}
   
   We next approximate $q_k(\lambda,N)$ by a Chebyshev polynomial of degree $v$, with the approximation error bounded by its total variation over $[-1,1]$, denoted by $V_{[-1,1]}$, according to Jackson's theorem (see e.g.,\cite{trefethen2013approximation,bernstein1912best}). 
   \begin{theorem}[\text{Chebyshev approximation on $[-1,1]$}, \cite{trefethen2013approximation}]
   \label{thm: Chebyshev}
	For a positive integer $k$, let $h : [-1,1] \to \mathbb{R}$ be a function such that its derivatives through $h^{(k-1)}$ are absolutely continuous on $[-1,1]$, and its $k$-th derivative $h^{(k)}$ is of bounded variation $V$. Then its Chebyshev  polynomial approximation $p_v$ of degree $v> k$ satisfies
	\begin{displaymath}
	  \|h-p_v\|_{[-1,1]} \leq \dfrac{4V}{\pi k(v-k)^k}.
	\end{displaymath}
	\end{theorem}
   
	\begin{proposition}\label{penalty}
	  For any integer $v >k$, there exists a non-negative polynomial  $p_k^{[v]}(\lambda,N)(t)$
	   of degree $v$  approximating the function
    $q_k(\lambda,N)(t)$ in \eqref{smooth penalty} from above  on $[-1,1]$ and satisfies that
	 \begin{displaymath}
	  	0\;\leq\; p_k^{[v]}(\lambda,N)(t)-q_k(\lambda,N)(t) \;\leq\; 8Ne^2a^{|\lambda|v} \quad \forall t \in [-1,1],
    \end{displaymath}
    where  $a= e^{-1/ (2e+1)}\approx 0.8561$.
	  \end{proposition}
	  \begin{proof}
           By Lemma~\ref{lemma: AC and BV}, $q_k(\lambda,N)$ satisfies the assumption of Theorem~\ref{thm: Chebyshev}. Then there exists a Chebyshev polynomial  $q_k^{[v]}(\lambda,N)$ which approximates $q_k(\lambda,N)$ with error satisfying
        	  	\begin{displaymath}
        	  	    \left\|q_k(\lambda,N)(t)-q_k^{[v]}(\lambda,N)(t)\right\|_{[-1,1]} \leq\; \frac{12N|\lambda|^{-k}4^kk!k}{\pi k(v-k)^k}\leq\; \frac{4N|\lambda|^{-k}4^kk!}{(v-k)^k}.
        	  	\end{displaymath}
                We next use the same idea as in Figure~\ref{fig:approximation} to vertically shift $q_k^{[v]}(\lambda,N)$ by the constant $\frac{4N|\lambda|^{-k}4^kk!}{(v-k)^k}$ to obtain a polynomial approximation from above as follows:
            \begin{displaymath}
	  	   p_k^{[v]}(\lambda,N)(t):= q_k^{[v]}(\lambda,N)(t)+\frac{4N|\lambda|^{-k}4^kk!}{(v-k)^k} \geq q_k(\lambda,N)(t) \geq 0 \quad \forall t \in [-1,1].
	  	\end{displaymath}
    Moreover, the approximation error on $[-1,1]$ is bounded by 
    \begin{displaymath}
	  	0\;\leq\; p_k^{[v]}(\lambda,N)(t)-q_k(\lambda,N)(t) \;\leq\; \frac{8N|\lambda|^{-k}4^kk!}{(v-k)^k} \quad \forall t \in [-1,1].
    \end{displaymath}
    According to Lemma~\ref{lem: choosing k}, we can choose $k <v$ in terms of $v$ and $\lambda$ such that 
    \begin{displaymath}
        \frac{|\lambda|^{-k}4^kk!}{(v-k)^k}\leq\; e^{2}a^{|\lambda| v}.
    \end{displaymath}
    From here, the required result follows.
    \end{proof}
    
   %%%%%%%%%%%%%%%%%%%%%%%%%
\subsection{A Hierarchy for Binary Polynomial Optimization Problems}\label{section: hierarchy for discrete}

We prove Theorem~\ref{thm: convergence rate for bpop} that gives a bound between $\lb(f,\widetilde{\cQ}(\cX))_r$ defined in \eqref{proposed hierarchy for bpop} and $f_{\min}$.

\begin{theorem}\label{thm: convergence rate for bpop}
Consider the problem \eqref{BPOP}. Assume that $\cX$ satisfies Assumption~\ref{normalization assumption} with $\cC = \bB^n$. For $r \in \mathbb{N}$, let $\xi_r^n$ be the least root of the degree-$r$ Krawtchouk polynomial with parameter $n$. If $(r+1)/n \leq 1/2$ and $d(d+1)\cdot \xi_{r+1}^n\leq 1/2$, then for any integer $v>k$ such that $l(v+1) \leq r$ and $(2l(v+1))(2l(v+1)+1)\xi_{r+1}^n \leq 1/2$, the following inequality holds:
    \begin{align*}
        0 \;&\leq\; f_{\min} -\lb(f,\widetilde{\cQ}(\cX))_r \\
        &\leq\; 2 \biggl(  \|f\|_{\bB^n}
        \Big(\dfrac{-2m}{\lambda}+1\Big)+ V(m,v,\lambda) \biggr)\cdot \max\{C_d,C_{2l(v+1)}\}\cdot \xi_{r+1}^n/n
        + V(m,v,\lambda), 
	\end{align*}
    where $V(m,v,\lambda) =\|f\|_{\bB^n} 16 \,m\,e^2   |\lambda|^{-1}a^{-\lambda v}$ with $a= e^{-1/ (2e+1)}\approx 0.8561$.
\end{theorem}
\begin{proof}
Let $(t, \sigma,h)$ be a feasible solution of \eqref{proposed hierarchy for bpop}, with $t \in \bR$, $\sigma \in \cQ(\bB^n)_{2r}$ and $h \in \cP_+([-1,1])_{\lfloor r/l \rfloor -1}$ satisfying that 
\begin{displaymath}
	f(x) - t \;=\; \sigma(x) + \langle h(G(x)),G(x) \rangle. 
\end{displaymath}
Then, for any $x \in \cX$, $\sigma(x) \geq 0$ and $\langle h(G(x)),G(x) \rangle=\sum_{j\in [m]}\lambda_j(x)h(\lambda_j(x)) \geq 0$. Hence, $t \leq f_{\min}$, and that implies $\lb(f,\widetilde{\cQ}(\cX))_r \leq f_{\min}$.

To get an upper bound on $f_{\min}-\lb(f,\widetilde{\cQ}(\cX))_r$, we proceed 
with the following scheme:
\begin{enumerate}
    \item Construct the penalty polynomial $P(x) = -\langle p_k^{[v]}(\lambda,N)(G(x)),G(x) \rangle$, with 
    $\lambda :=\max_{\bB^n\backslash\cX}\lambda_m(x)$, 
    $N= -2\|f\|_{\bB^n}/\lambda$,
    and $k$ as in Section~\ref{lem: choosing k}. Consider the penalized POP:
    \begin{equation}\label{bpop1}
        F_{\min}:=\;\min_{x \in \bB^n}F(x) := f(x) +P(x).
    \end{equation}
    \item We define 
    \begin{equation}\label{hierarchy: bpop1 F}
        \lb(F,\cQ(\bB^n))_r :=\; \max\{F-t:\; F-t \in \cQ(\bB^n)_{2r}\},
    \end{equation}
    and prove that $\lb(f,\widetilde{\cQ}(\cX))_r \geq \lb(F,\cQ(\bB^n))_r$. Then we bound
    \begin{equation}\label{eq: bpop1 1}
        f_{\min}-\lb(f,\widetilde{\cQ}(\cX))_r \leq\; f_{\min}-F_{\min} +F_{\min}-\lb(F,\cQ(\bB^n))_r.
    \end{equation}
    \item We bound $f_{\min}-F_{\min}$ based on the construction of $p_k^{[v]}$.
    \item We bound $F_{\min}-\lb(F,\cQ(\bB^n))_r$ by Theorem~\ref{binary_rate}.
\end{enumerate}
Now we elaborate each step in the above scheme. 

\textbf{Step 1}: Under Assumption~\ref{normalization assumption} with $\cC = \bB^n$, we know $\lambda \in [-1,0)$. 
Consider the penalty function $q_k(\lambda,N)$ as in \eqref{smooth penalty} in 
Section~\ref{section: Penalty for binary},
we have that $p_k^{(v)}(\lambda,N)$ is well-defined as in Proposition~\ref{penalty}. Hence 
$P(x)$ is also well-defined with degree bounded by $2l(v+1) \leq 2r$, since $l(v+1) \leq r$.

\textbf{Step 2}: We prove the inequality $\lb(f,\widetilde{\cQ}(\cX))_r \geq \lb(F,\cQ(\bB^n))_r$. Let $(t,\sigma)$ be a feasible solution of \eqref{hierarchy: bpop1 F}, that is
\begin{displaymath}
    F(x) -t = \sigma(x) \in \cQ(\bB^n)_{2r} \; \Rightarrow \; f(x) -t =\; \sigma(x) +\langle p_k^{[v]}(\lambda,N)(G(x)),G(x) \rangle \in \widetilde{\cQ}(\cX)_{2r}.
\end{displaymath}
Therefore, $(t,\sigma, p_k^{[v]}(\lambda,N))$ is a feasible solution of \eqref{proposed hierarchy for bpop}, which implies that 
\begin{displaymath}
  \lb(f,\widetilde{\cQ}(\cX))_r \geq t \; \Rightarrow \;  \lb(f,\widetilde{\cQ}(\cX))_r \geq \lb(F,\cQ(\bB^n))_r. 
\end{displaymath}
\textbf{Step 3}: For evaluation convenience, we define the function 
\begin{equation}\label{eq: bpop1 2}
    \left\langle q_k\left(\lambda,N\right)(G(x)),G(x) \right\rangle :=\; \sum_{j \in [m]}\lambda_j(x)q_k\left(\lambda,N\right)(\lambda_j(x)).
\end{equation}
Notice that we slightly abuse the notation here since $q_k\left(\lambda,N\right)$ is not a polynomial. Since $p_k^{[v]}(\lambda,N)$ is a polynomial approximation of $q_k(\lambda,N)$, \eqref{eq: bpop1 2} serves as a tool to handle the evaluations in our proof. Consider the following function
\begin{displaymath}
    \overline{F}(x)= f(x)-\left\langle q_k\left(\lambda,N\right)(G(x)),G(x)\right\rangle.
\end{displaymath}
We first estimate the gap between $F$ and $\overline{F}$ over $\bB^n$. Using Lemma~\ref{penalty}, we get for any $x \in \bB^n$,
\begin{eqnarray}\label{binary rate inequality 1}
     \left|\overline{F}(x) - F(x)\right| &\leq&  \sum_{i=1}^m \left|\left(p_k^{[v]}\left(\lambda,N\right)(\lambda_i(x))-q_k\left(\lambda,N\right)(\lambda_i(x))\right)\lambda_i(x)\right|
     \nonumber \\
    & \leq& m\cdot 8 Ne^2a^{|\lambda|v} =  \|f\|_{\bB^n} 16 \,m\,e^2   |\lambda|^{-1}a^{-\lambda v} = V(m,v,\lambda).
\end{eqnarray}
From \eqref{binary rate inequality 1}, we get
\begin{eqnarray}\label{binary rate inequality 2}
&\|F-\overline{F}\|_{\bB^n} & \;\leq\;  V(m,v,\lambda).
\end{eqnarray}
We next show that $\min_{x\in\mathbb{B}^n}\overline{F}(x) =\min_{x\in\mathcal{X}}\overline{F}(x)= f_{\min}$. 
Recall the definition of $q_k\left(\lambda,N\right)$  in \eqref{smooth penalty}. For any $x \in \mathcal{X}$, all the eigenvalues $\lambda_i(x)$ are non-negative and at most $1$ by Assumption~\ref{normalization assumption}. Thus we have 
\begin{displaymath}
    \left\langle q_k\left(\lambda,N\right)(G(x)),G(x) \right\rangle= \sum_{i=1}^m q_k\left(\lambda,N\right)(\lambda_i(x))\lambda_i(x) =0\quad \forall x \in \cX.
\end{displaymath}
    Consequently, $f(x) = \overline{F}(x)\ \forall x \in \mathcal{X}$, which implies 
    that $\min_{x\in\mathcal{X}}\overline{F}(x) = f_{\min}.$ 
    Now for any point $x \in \mathbb{B}^n \backslash \mathcal{X}$, $\lambda_m(x) \in [-1,\lambda]$. Whence, we get
\begin{eqnarray*}
    \overline{F}(x) &=& f(x)-\left\langle q_k\left(\lambda,N\right)(G(x)),G(x)\right\rangle \\
    &=&\ f(x) -q_k\left(\lambda,N\right)(\lambda_m(x))\lambda_m(x)-\sum_{i=1}^{m-1} q_k\left(\lambda,N\right)(\lambda_i(x))\lambda_i(x)\\
    &\geq &\ f(x) - N \cdot \lambda = f(x) +2\|f\|_{\bB^n} \geq \|f\|_{\bB^n} \geq f_{\min}.
\end{eqnarray*}
The first inequality is due to the fact that for any $i \in [m-1]$, $q_k\left(\lambda,N\right)(\lambda_i(x))\lambda_i(x) =0 $ if $\lambda_i(x)\geq 0$, and $q_k\left(\lambda,N\right)(\lambda_i(x))\lambda_i(x) \leq 0 $ if $\lambda_i(x)\leq 0$. Therefore, for any $x \in \mathbb{B}^n \backslash\mathcal{X}$, we obtain that $\overline{F}(x) \geq f_{\min}$. Whence, $\min_{x\in\mathbb{B}^n}\overline{F}(x) = \min_{x\in\mathcal{X}}\overline{F}(x) = f_{\min}$.

Next, we show that $|F_{\min}-f_{\min}|\leq V(m,v,\lambda)$.
Let $\bar{x}\in\cX$ be a minimizer such that 
$\overline{F}(\bar{x}) = \min_{x\in\cX}\overline{F}(x) = f_{\min}$.
First, we have $F_{\min} - f_{\min}\leq F(\bar{x})-\overline{F}(\bar{x}) \leq V(m,v,\lambda)$, by \eqref{binary rate inequality 1}.
Then, letting $x^*$ be a minimizer of $F$ over $\mathbb{B}^n$, i.e., $F(x^*)=F_{\min}$, we have 
\begin{equation*}
   F_{\min}-f_{\min}
   = F(x^*) -\overline{F}(x^*) + \overline{F}(x^*)-f_{\min}
   \geq F(x^*) -\overline{F}(x^*) 
   \geq - V(m,v,\lambda),
\end{equation*}
where the first inequality follows from the fact that $\overline{F}(x^*)\geq f_{\min}$
and the second inequality follows from \eqref{binary rate inequality 1}.
Thus we have shown that 
\begin{eqnarray}
    |F_{\min}-f_{\min}|\leq V(m,v,\lambda).
    \label{eq: bpop fmin to Fmin}
\end{eqnarray}

\textbf{Step 4:} Here we apply Theorem \ref{binary_rate} to bound $F_{\min}- \lb(F,\cQ(\bB^n)_r$. Before doing so, we need to check the conditions of Theorem \ref{binary_rate} on $F$. We begin with the degree of $F$. Since the degree of $G(x)$ is bounded by $2l$, the degree of $\langle p_k^{[v]}(\lambda,N)(G(x)),G(x) \rangle$ is bounded by $2l(v+1)$. Then, $\deg(F) \leq \max\{d,2l(v+1)\}$. The inequalities  
\begin{displaymath}
   d(d+1)\xi_{r+1}^n \leq 1/2,\quad \mbox{and} \quad (2l(v+1))(2l(v+1)+1)\xi_{r+1}^n \leq 1/2
\end{displaymath}
implies that $\deg(F)(\deg(F)+1)\xi_{r+1}^n \leq 1/2$, which satisfies the degree condition of Theorem~\ref{binary_rate}. We also need to estimate $\|F\|_{\mathbb{B}^n}$. For all $x \in \mathbb{B}^n$, since $\lambda_i(x) \in [-1,1]$ for all $i \in [m]$,  we can upper bound $\|F\|_{\bB^n}$ as follows: for all $x \in \mathbb{B}^n$, we have
\begin{equation*}
    \left|\overline{F}(x)\right| = \left|f(x) -\sum_{i=1}^{m}
    q_k\left(\lambda,N\right)(\lambda_i(x))\lambda_i(x)\right| \leq \|f\|_{\bB^n} + mN = \|f\|_{\bB^n}\left(1 - \dfrac{2m}{\lambda}\right).
\end{equation*}
Then we recall the inequality \eqref{binary rate inequality 2} to obtain the bound 
\begin{equation*}
\|F\|_{\bB^n} \leq \|\overline{F}(x)\|_{\bB^n} + \|F(x) -\overline{F}(x)\|_{\bB^n} \leq \|f\|_{\bB^n}\left(1-\dfrac{2m}{\lambda} \right)+ V(m,v,\lambda).
           %%+\dfrac{16m \cdot M_f\cdot M_k \cdot |\lambda|^{-k-1}}{\pi k(v-k)^k}.
\end{equation*}
Now we can apply Theorem~\ref{binary_rate} to obtain the inequality 
\begin{equation}\label{eq: bpop Fmin and lbF}
    F_{\min} - \lb(F,\mathcal{Q}(\mathbb{B}^n))_r \leq  2\left[\|f\|_{\bB^n}\left(1-\dfrac{2m}{\lambda}\right) + V(m,v,\lambda)
           %+\dfrac{16m \cdot M_f\cdot M_k \cdot |\lambda|^{-k-1}}{\pi k(v-k)^k}
\right]\max\{C_d,C_{vl+l}\}\cdot \xi_{r+1}^n/n.
\end{equation}
Finally, we substitute the bounds of \eqref{eq: bpop fmin to Fmin} and \eqref{eq: bpop Fmin and lbF} into \eqref{eq: bpop1 1} to complete the proof.
\end{proof}
\begin{remark}
    %While Corollary~\ref{thm: convergence for bpop} ensures the convergence of %\eqref{proposed hierarchy for bpop}, 
    Theorem~\ref{thm: convergence rate for bpop} provides a bound when $2r \leq n$. For $2r > n$, $\cQ(\bB^n)= \cQ(\bB^n)_{2r}$, and hence $F_{\min}= \lb(F,\cQ(\bB^n))_r$. Therefore, the bound for the regime $2r \geq n$ is reduced to 
    \begin{displaymath}
        f_{\min} - \lb(f,\widetilde{\cQ}(\cX))_r \leq 
        \|f\|_{B^n} 16m e^2|\lambda|^{-1}a^{-\lambda v}.
    \end{displaymath}
    We can choose $v = \lfloor r/l \rfloor -1$ so that $\lb(f,\widetilde{\cQ}(\cX))_r$ converges to $f_{\min}$ at the exponential rate of $\mathcal{O}(a^{-\lambda r/l})$.
\end{remark}

   %%%%%%%%%%%%%%%%%%%%%%%%%%
\section{Continuous Polynomial Optimization Problems with semidefinite constraints}\label{section: Continuous}

By scaling $G(x)$ if necessary, we may assume that the  domain $\mathcal{X}$ in \eqref{SPOP} is contained in the unit ball $B^n$. In this section, we consider the  following problem:
\begin{align}\label{cpop}
    f_{\min} = \min{f(x)}\quad \mbox{subject to } x \in \mathcal{X} := \{ x\in B^n \,:\, G(x)\succeq 0\}. 
\end{align}
Under Assumption~\ref{normalization assumption} with $\cC = B^n$, we propose the hierarchy \eqref{proposed hierarchy for cpop}, whose corresponding 
dual SDP problem  is given in the following corollary.
\begin{corollary}\label{cor: proposed momemnt SDP for cpop}
The dual of the SDP problem corresponding to \eqref{proposed hierarchy for cpop} is given by
\begin{align}\label{proposed momemnt SDP for cpop}
    \min_{y\in \bR^{s(n,2r)}}&\quad L_y(f)\\
    \mbox{subject to}&\quad y_0 =1,\; \M_r(y) \succeq 0,\; \M_{r-1}\biggl(\bigl(1-\sum_{i=1}^nx_i^2\bigr)y\biggr)\succeq 0,\nonumber\\
    &\quad L_y(P_{\lceil(\lfloor r/l \rfloor -1)/2 \rceil}^G)\succeq 0,\; L_y(Q_{\lceil(\lfloor r/l \rfloor -1)/2 \rceil-1}^G) \succeq 0.\nonumber
\end{align}
\end{corollary}
\begin{proof}
    The proof is analogous to the proof of Corollary~\ref{cor: proposed momemnt SDP for bpop} via using Proposition~\ref{prop: SDP certificate}.
\end{proof}
\begin{remark}
In contrast to the hierarchy of Hol and Scherer 
in \eqref{hierarchy: Hol and Scherer},
    we can observe that the hierarchy \eqref{proposed momemnt SDP for cpop} avoids the explosion of the size of the SDP matrix.
\end{remark}

In what follows, we study the convergence rate of \eqref{proposed hierarchy for cpop} by the same penalty function framework proposed in Section~\ref{section: hierarchy for discrete} with several modifications. Since we do not have the parameter $\lambda$ as in the discrete setting 
for \eqref{cpop}, we consider a sequence of small neighborhoods $\{V(\delta)\}_{\delta >0}$ of $\mathcal{X}$ such that $\max_{x \in V(\delta)}\dist(x,\mathcal{X}) \to 0$ as $\delta \to 0$. For each $\delta >0$, the penalty polynomial $P$ will lift up the value of  $f$ outside $V_{\delta}$, but slightly perturb the value of $f$ over $\mathcal{X}$. The "middle area", defined as $V_{\delta}\backslash\mathcal{X}$, allows us to modify the penalty function as smoothly as possible. According to Jackson's theorem, the smoother a function is, the smaller is the polynomial approximation error. However, the degree of $P$ will tend to infinity as $\delta \to 0$. To resolve this obstacle, we leverage the Lipschitz constant of $f$ and the Łojasiewicz exponent to establish the connection between $\delta$ and the degree of the penalty function. 
To construct the neighborhood $V(\delta)$ of $\mathcal{X}$, 
we need to make use of the 
Łojasiewicz inequality for a set defined by a matrix inequality as stated in 
Theorem~\ref{Lojasiewicz}.

\begin{comment}
We review a version of the Łojasiewicz inequality for matrix inequality (see e.g., \cite[Theorem 4.1]{dinh2016lojasiewicz}).

\begin{theorem}\label{Lojasiewicz-2}
For any compact set $K$ containing $\cX$, there exist a Łojasiewicz constant $C>0$ and a Łojasiewicz exponent $0 < L \leq 1$ depending on $\mathcal{X}$ such that
\begin{displaymath}
    d(x,\mathcal{X})\leq C\max\{0,-\lambda_m(x)\}^{L}\quad \forall x \in K.
\end{displaymath}
\end{theorem}
\begin{remark}
Note that the condition $L \leq 1$ is not present in \cite{dinh2016lojasiewicz}. However, we can assume it to be at most $1$ because of the compactness of $K$. In particular, if $L >1$, we can replace $L$ by $1$ and multiply the Łojasiewicz constant by $\max_{x \in K}\max\{0,-\lambda_m(x)\}^{L-1} <\infty$ to obtain a new inequality with the Łojasiewicz exponent $1$.
\end{remark}
\end{comment}

For any $\delta \in (0,1]$, we define the neighborhood $V(\delta)$ of $\mathcal{X}$ by
\begin{displaymath}
V(\delta) \;=\; \{x \in B^n: G(x)+\delta I_m \succeq 0\} \;=\; \{x \in B^n: \lambda_m(x) \geq -\delta \}.
\end{displaymath}
We apply Theorem~\ref{Lojasiewicz} to bound the distance of a point in $V(\delta)$ to $\mathcal{X}$ as follows:
\begin{eqnarray} \label{Loj-ineq}
    \max_{x \in V(\delta)}\dist(x,\mathcal{X}) \leq C\cdot \delta^L.
\end{eqnarray}
We next state and prove the bound between $f_{\min}$ and $\lb(f,\overline{\cQ}(\cX))_r$.
\begin{theorem}\label{convergence_thm_cpop}
	%Let $f$ be a polynomial of degree $d$, $l$ is the degree of $G(x)$. The domain $\mathcal{X}$ is defined by a matrix inequality as in~\ref{cpop}. The optimal value of $f$ over $\mathcal{X}$ is denoted by $f_{\min}$. 
 Consider the problem \eqref{cpop} with $\cX$ contained in $B^n$.
 Let $L_f$ be the Lipschitz number of $f$ over $B^n$, and $C$ and $L$ are the Łojasiewicz constant and exponent for $\mathcal{X}$ as in Theorem~\ref{Lojasiewicz} with $K=B^n$.
 Let $\delta$ be an arbitrary positive number in $(0,1]$. 
 For any positive integer $r$ and $v$ such that $r \geq 2n\max\{d,2l(v+1)\}$, $\lb(f,\overline{\cQ}(\cX))_r$ is a lower bound of $f_{\min}$ satisfying that
\begin{eqnarray}\label{final inequality cpop}
    f_{\min} - \lb(f,\overline{\cQ}(\cX))_r
    &\leq& 2\Big( \|f\|_{B^n}\big(1 + \frac{m}{\delta} \big)+U(m,v,\delta)\Big)
    \cdot\frac{C_B(n,\max\{d,2l(v+1)\})}{r^2}
    \nonumber\\
   && +\; U(m,v,\delta)+ CL_f\delta^L, 
\end{eqnarray}
where $U(m,v,\delta) = \|f\|_{B^n} 16e^2m\delta^{-1} a^{\delta v}.$
\end{theorem}
	  
\begin{proof} This proof uses several similar arguments 
as in the proof of Theorem~\ref{thm: convergence rate for bpop}. For instance, 
by replacing $\mathbb{B}^n$ by $B^n$ in the early part of the 
proof of Theorem~\ref{thm: convergence rate for bpop}, 
we have
 \begin{displaymath}
     \lb(f,\overline{\cQ}(\cX))_r \leq f_{\min}.
 \end{displaymath}
To prove the upper bound on $f_{\min}-\lb(f,\overline{\cQ}(\cX))_r$, we use 
the following scheme:
\begin{enumerate}
    \item Construct the penalty polynomial $P(x) = -\langle p_k^{[v]}(-\delta,N)(G(x)),G(x) \rangle$ for a fixed $\delta >0$ and $N= 2\|f\|_{B^n}/\delta$. Consider the pernalized POP:
    \begin{equation}\label{cpop1}
        F_{\min}:=\;\min_{x \in B^n}F(x)\;:= f(x) +P(x).
    \end{equation}
    \item We define 
    \begin{equation}\label{hierarchy: cpop1 F}
        \lb(F,\cQ(B^n))_r :=\; \max\{F-t:\; F-t \in \cQ(B^n)_{2r}\},
    \end{equation}
    and prove that $\lb(f,\overline{\cQ}(\cX))_r \geq \lb(F,\cQ(B^n))_r$. Then we bound
    \begin{equation}\label{eq: cpop 1}
        f_{\min}-\lb(f,\overline{\cQ}(\cX))_r \leq\; f_{\min}-F_{\min} +F_{\min}-\lb(F,\cQ(B^n))_r.
    \end{equation}
    \item We bound $f_{\min}-F_{\min}$ based on the construction of $p_k^{[v]}$.
    \item We bound $F_{\min}-\lb(F,\cQ(B^n))_r$ by Theorem~\ref{unit ball}.
\end{enumerate}
Although Step $1$ and $2$ can be proved similarly to the proof of Theorem~\ref{thm: convergence rate for bpop}, Step $3$ and $4$ require more complicated arguments, which we provide next. 

\textbf{Step 3}: We define the following function as a tool for our proof:
	  	\begin{displaymath}
	  	    \overline{F}(x)= f(x)-\left\langle q_k\left(-\delta,N\right)(G(x)),G(x) \right\rangle,
	  	\end{displaymath}
where $\left\langle q_k\left(-\delta,N\right)(G(x)),G(x) \right\rangle$ is defined as in \eqref{eq: bpop1 2}. Note that $q_k\left(-\delta,N\right)$ has the following property: 
\begin{eqnarray} \label{eq-q}
q_k\left(-\delta,N\right)(t) \left\{
\begin{array}{ll}
 = 0 & \mbox{if $t \in [0,1]$}
\\
\geq 0 & \mbox{if $t \in [-1,0]$.}
\end{array}\right.
\end{eqnarray}
Based on \eqref{eq-q}, we have that $\overline{F}(x)\geq f(x)$ for all $x\in B^n.$

Under Assumption~\ref{normalization assumption} with $C = B^n$, applying Lemma~\ref{penalty} with $N= 2\|f\|_{B^n}/\delta$ and $\lambda = -\delta$ implies that 
\begin{equation}\label{cpop inequality 1}
    |F(x) -  \overline{F}(x)| \leq 8e^2mNa^{|\lambda|v}=16e^2m\|f\|_{B^n}\delta^{-1} a^{\delta v}=U(m,v,\delta) \quad \forall\; x\in B^n.
\end{equation}     
This induces the gap between the minimum values of $F$ and $\overline{F}$ over $B^n$ as
\begin{equation}\label{eq-gapFr}
\Big|F_{\min} -  \min_{x\in B^n}\overline{F}(x)\Big| \;\leq\; U(m,v,\delta).
\end{equation}
We next estimate the gap between  $\min_{x\in B^n}\overline{F}$ and $f_{\min}$
by considering three cases. 
\\
(i) For any $x \in \mathcal{X}$, all the eigenvalues of $G(x)$ are contained in $[0,1]$, which leads to
\begin{displaymath}
      \left\langle q_k\left(-\delta,N\right)(G(x)),G(x) \right\rangle= \sum_{i=1}^mq_k\left(-\delta,N\right)(\lambda_i(x))\lambda_i(x)=0.
\end{displaymath}
This implies that $\overline{F}(x) = f(x)\ \forall x \in \mathcal{X}$.

\noindent(ii) For $x \in B^n \backslash V(\delta)$, we have that $\lambda_m(x) \in [-1,-\delta]$.  Using the property of $q_k(-\delta,N)$ in \eqref{eq-q}, we obtain that 
       \begin{align*}
           &\overline{F}(x) = f(x)-\left\langle q_k\left(-\delta,N\right)(G(x)),G(x)\right\rangle \\
           =\;& f(x) -q_k\left(-\delta,N\right)(\lambda_m(x))\lambda_m(x)-\sum_{i=1}^{m-1} q_k\left(-\delta,N\right)(\lambda_i(x))\lambda_i(x)\\
           \geq\; & f(x) + \delta\cdot N = f(x) +2\|f\|_{B^n}\geq  \|f\|_{B^n} \geq f_{\min}.
       \end{align*}
       (iii) For any $x \in V(\delta)\backslash\mathcal{X}$, we choose $\overline{x} \in \mathcal{X}$ to be a closest point to $x$ in $\mathcal{X}$. Then, the Łojasiewicz inequality 
       \eqref{Lojasiewicz} implies that
       \begin{equation*}\label{eq: min of F bar}
            \overline{F}(x) - f_{\min} \geq f(x) -f_{\min} \geq f(x) - f(\overline{x}) \geq -L_f\dist(x,\overline{x}) \geq -L_f C\delta^L.
       \end{equation*}
        Combining all the three cases, we obtain that 
       \begin{displaymath}
           f_{\min} - \min_{x\in B^n}\overline{F}(x)  \leq L_f C\delta^L.    
       \end{displaymath}
   Together with \eqref{eq-gapFr}, the gap between $F_{\min}$ and $f_{\min}$ satisfies that
        \begin{equation}\label{cpop inequality 2}
            f_{\min}-F_{\min} =  f_{\min}- \min_{x\in B^n}\overline{F}(x) + 
            \min_{x\in B^n}\overline{F}(x) - F_{\min}
            \;\leq\;
            L_f C\delta^L + U(m,v,\delta).
            %%\dfrac{16m \cdot M_f \cdot M_k \cdot \delta ^{-k-1}}{\pi(v-k)^k}+ L_f C\delta^L.
        \end{equation}
        
\textbf{Step 4}: We next bound $F_{\min}-\lb(F,\cQ(B^n))_r$ by applying Theorem~\ref{unit ball} to \eqref{hierarchy: cpop1 F}. Before doing so, we need to check the condition of Theorem~\ref{unit ball} for $F(x)$ over $B^n$. Similar to the proof of Theorem~\ref{thm: convergence rate for bpop}, the degree of $F$ is upper bounded by $\max\{d,2l(v+1\}$. The condition $r \geq 2n\max\{d,2l(v+1)\}$ implies that $r \geq 2n \deg (F)$. 

We estimate the maximum value of $F(x)$ on $B^n$ as follows: by the definition of $q_k(-\delta,N)$, for all $x \in B^n$, we have
\begin{displaymath}
\left|\left\langle q_k\left(-\delta,N\right)(G(x)),G(x) \right\rangle\right|\leq  \sum_{i=1}^m\left|q_k\left(-\delta,N\right)(\lambda_i(x))\lambda_i(x)\right| \leq mN = \dfrac{2m\|f\|_{B^n}}{\delta}.
\end{displaymath}
        Using the inequality in \eqref{cpop inequality 1} and the above inequality,
        we get
        \begin{displaymath}
            \max_{x\in B^n}F(x) \leq \max_{x\in B^n}f(x) + \dfrac{2m\|f\|_{B^n}}{\delta}
            + U(m,v,\delta).
            %%+\dfrac{16m \cdot M_f \cdot M_k \cdot \delta ^{-k-1}}{\pi(v-k)^k}.
        \end{displaymath}
         Using the inequality in \eqref{cpop inequality 1} again, we have
        \begin{align*}
        U(m,v,\delta)  + F_{\min} &\geq \min_{x\in B^n} 
        \Big(\overline{F}(x) = f(x) - \sum_{i=1}^{m} q_k\left(-\delta,N\right)(\lambda_i(x))\lambda_i(x) \Big)
        %\\
        %& 
        \geq\;  \min_{x\in B^n} f(x), 
        \end{align*}
        where the last inequality used the fact that 
        $\sum_{i=1}^{m} q_k\left(-\delta,N\right)(\lambda_i(x))\lambda_i(x) \leq 0$ 
        for all $x\in B^n$.
        Thus, we have a bound on the gap between the maximum and  minimum values of $F$ over $B^n$ as 
        \begin{eqnarray*}\label{eq: max of F}
            \max_{x\in B^n}F(x) - F_{\min} 
            &\leq & \max_{x\in B^n}f(x) - \min_{x\in B^n}f(x) + \dfrac{2m\|f\|_{B^n}}{\delta}+ 2 U(m,v,\delta)
            \nonumber \\
           &\leq& 2 \|f\|_{B^n} \Big(1 +
            \dfrac{m}{\delta}\Big)+ 2 U(m,v,\delta).
        \end{eqnarray*}
    Hence, applying Theorem~\ref{unit ball} to $F$ provides the following bound:
    \begin{align}\label{eq: F and lb cpop}
        0 &\leq F_{\min}-\lb(F,\mathcal{Q}(B^n))_r \nonumber\\
        \leq & \; 2\left(\|f\|_{B^n} \Big(1 +
            \dfrac{m}{\delta}\Big)+  U(m,v,\delta)\right)\dfrac{C_B(n,\max\{d,2l(v+1)\})}{r^2}.
    \end{align}
    Using the bounds of \eqref{eq: F and lb cpop} and \eqref{cpop inequality 2} in \eqref{eq: cpop 1}, we get our desired inequality \eqref{final inequality cpop}.
\end{proof}

\begin{remark}
   We can replace the containment of the set $\{x\in \mathbb{R}^n\mid G(x) \succeq 0\}$ in the unit ball \eqref{cpop} 
	  by containment in other simple sets, such as the hypercube and standard simplex, which exhibit the same convergence rate for the SOS hierarchy of lower bounds as that for the unit ball. Consequently, different hierarchies of lower bounds can be constructed with convergence rates matching that in Theorem~\ref{convergence_thm_cpop}. However, we prefer the unit ball for its simplicity — it requires only a single polynomial of degree $2$ for its description, which enhances the computational efficiency of our proposed SOS hierarchy.
	  \end{remark}
\begin{remark}
In the case when $G$ has a block-diagonal structure, that is, 
$G= \diag (G_1,\dots,G_k)$ for matrices $G_i \in \bS\bR[x]^{m_i}$ (we let $l_i = \lceil \deg(G_i) /2 \rceil$), the feasible set $\cX$ is defined by 
    \begin{displaymath}
        \mathcal{X}= \{ x \in B^n: \ G_i(x) \succeq 0,\ i \in [k]\}.
    \end{displaymath}
We can either apply the hierarchy \eqref{proposed hierarchy for cpop} to approximate \eqref{cpop}, or a new hierarchy based on the block-diagonal structure as follows: 
\begin{eqnarray}\label{proposed hierarchy for multiple matrices}
    &\lb(f,\overline{\cQ}(G_1,\ldots,G_k))_r =\; \max\{ t:\; f-t \in \overline{\cQ}(G_1,\ldots,G_k)_{2r}\},
    \\
    & \overline{\cQ}(G_1,\ldots,G_k)_{2r} =\; \cQ(B^n)_{2r}+ \cH(G_1)_{\lfloor r/l_1 \rfloor-1}+\ldots +\cH(G_k)_{\lfloor r/l_k \rfloor-1}.\nonumber
\end{eqnarray}
We observe that this hierarchy is stronger than \eqref{proposed hierarchy for cpop} on $\cX$ since 
\begin{displaymath}
    \overline{\cQ}(\cX)_{2r}= \overline{\cQ}(G)_{2r}\subset \; \overline{\cQ}(G_1,\ldots,G_k)_{2r}.
\end{displaymath}
Moreover, the SDP relaxation corresponding to \eqref{proposed hierarchy for multiple matrices} is given by
\begin{align}\label{proposed momemnt SDP for multiple matrices}
    \min_{y\in \bR^{s(n,2r)}}&\quad L_y(f)\\
    \mbox{subject to}&\quad y_0 =1,\; \M_r(y) \succeq 0,\; \M_{r-1}\biggl(\bigl(1-\sum_{i=1}^nx_i^2\bigr)y\biggr)\succeq 0,\nonumber\\
    &\quad L_y(P_{\lceil(\lfloor r/l_i \rfloor -1)/2 \rceil}^{G_i})\succeq 0,\quad L_y(Q_{\lceil(\lfloor r/l_i \rfloor -1)/2 \rceil-1}^{G_i}) \succeq 0 \;\; \forall\;\; i \in [k],\nonumber
\end{align}
which also avoids the explosion of the size of the SDP matrices. Since \eqref{proposed hierarchy for multiple matrices} is stronger than \eqref{proposed hierarchy for cpop}, the bound in Theorem~\ref{convergence_thm_cpop} is valid for $\lb(f,\overline{\cQ}(G_1,\ldots,G_k))_r$ under  Assumption~\ref{normalization assumption} with $\cC = B^n$. 

For the special case where $m_i =1 $ for all $i \in [k]$, $\eqref{cpop}$ is a scalar POP. Therefore, we have a new SDP relaxation 
\eqref{proposed momemnt SDP for multiple matrices} for a scalar POP.
%exploring the structure of univariate polynomials.
\end{remark}

%%===================================
\section{A version of the Putinar's Positivstellensatz}\label{section: another version}
This section focuses on using an identical framework as in the last section to develop another version of Putinar's Positivstellensatz for semialgebraic sets defined by semidefinite matrix constraints. This effort aims to establish analogous polynomial bounds on the degree as those in Theorem~\ref{thm-Baldi} and Theorem~\ref{thm-Baldi2} of the papers \cite{baldi2021moment,baldi2025lojasiewicz}. The main theorem of this section is  stated as follows.
    \begin{theorem}\label{thm: bound on degree}
    Assume without loss of generality that $\|f\|_{B^n} =1$.
    Let $L_f$ be the Lipschitz number of the degree-$d$ polynomial $f$ over $B^n$, and $C$ and $L$ are the Łojasiewicz constant and exponent for $\mathcal{X}$ as in Theorem~\ref{Lojasiewicz} with $K=B^n$. Set $d' = \max\{ \lceil d/2\rceil, l\}$.
    \\
(i)  Under the same setting of Theorem~\ref{convergence_thm_cpop}, 
the hierarchy $\lb(f,\overline{\cQ}(\cX))_r$ in \eqref{proposed hierarchy for cpop}
converges to $f_{\min}$ with the polynomial rate of $\mathcal{O}(r^{\frac{-2L}{n+2L+5}})$. 
\\
(ii) Suppose $f$ is positive over $\cX$ so that $E_f := f_{\min} > 0$. Let $\overline{r}$ be the smallest integer such that $\bar{r} \geq 2n\max\{d,2l(v+1)\}$.
Define $\delta \in (0,1]$ and $v$ by
         \begin{align*}
             &b = \frac{-2}{\ln (a)}\ln(r), \quad \gamma = \frac{(n+1)^2e^{\frac{n-2}{2}}d'^{\frac{n+3}{2}}}{CL_f},\\
             &v = \left\lfloor\gamma^{\frac{-2}{n+2L+5}} b^{\frac{2(L+1)}{n+2L+5}}r^{\frac{4}{n+2L+5}}\right \rfloor,\ \delta =\gamma^{\frac{2}{n+2L+5}} b^{\frac{n+3}{n+2L+5}}r^{\frac{-4}{n+2L+5}}.
         \end{align*}
         Then there exists a positive constant $\rho$ (independent of $n,m,l,L_f, C,L$) such that  $f \in \overline{\cQ}(\cX)_{2r}$ for any $r$ satisfying 
         \begin{displaymath}
             r \geq \max\{\overline{r},\ \big(\rho mn\max\{CL_f,1\}\big)^{1+\frac{n+5}{2L}}(ed')^{\frac{n+2L+5}{2}}E_f^{-1-\frac{n+5}{2L}}\}.
         \end{displaymath}
    
	\end{theorem}

    \begin{proof}
     (i) We first prove the convergence rate of the hierarchy of lower bounds \eqref{proposed hierarchy for cpop} by simplifying the right-hand side of \eqref{final inequality cpop}, which is denoted as 
    \begin{eqnarray}\label{bound cpop}
        W(v,\delta,r) &:=& 
        2\left( 1 + m\delta^{-1}+16e^2m\delta^{-1} a^{\delta v}\right)\dfrac{C_B(n,\max\{d,2l(v+1)\})}{r^2}
        \nonumber \\
        && +16e^2m\delta^{-1} a^{\delta v}+ CL_f\delta^L.
    \end{eqnarray}
     We observe that in order for the term in \eqref{bound cpop} to tend to $0$, $\delta v$ must tend to $\infty$ and $\delta$ must tend to $0$, which means that $v$ must tend to $\infty$. Under this scenario, we make the following simplification on the behavior of each of the terms in \eqref{bound cpop} using the the big-O notation $\mathcal{O}$ as follows:
    \begin{align*}
         & 2\left( 1 + m\delta^{-1}+16e^2m\delta^{-1} a^{\delta v}\right) = \mathcal{O}(m\delta^{-1})\quad \text{(since $a < 1$ and $\delta v \to \infty$)}\nonumber\\
        &C_B(n,\max\{d,2l(v+1)\}) = \mathcal{O}\big((n+1)^2e^{\frac{n-2}{2}}d'^{\frac{n+3}{2}}v^{\frac{n+3}{2}}\big)
        \quad  \text{(by Lemma~\ref{lem: parameter C})}\nonumber\\
        &16e^2m\delta^{-1} a^{\delta v} = \mathcal{O}(m\delta^{-1}a^{\delta v}).
    \end{align*}
    Hence
    \begin{equation} \label{approximation on bound}
        W(v,\delta,r)   = \mathcal{O} \big(m(n+1)^2e^{\frac{n-2}{2}}d'^{\frac{n+3}{2}}\delta^{-1}v^{\frac{n+3}{2}}r^{-2}+m\delta^{-1}a^{\delta v}+CL_f\delta^L\big).  
    \end{equation}
    We note that $\delta^{-1}v^{\frac{n+3}{2}}r^{-2}$ tends to $0$ slower than $r^{-2}$ since  $v >1$, and $\delta^{-1}a^{\delta v}$ tends to $0$ exponentially fast since
    $\delta v \to \infty$ and $a<1$. We can simplify the above bound by setting 
    \begin{eqnarray*}\label{delta equation 1}
       & &a^{\delta v} = r^{-2} \;\Leftrightarrow\; \delta  = bv^{-1}, \quad \text{where} \quad b := \frac{-2}{\ln (a)}\ln(r) \geq 0; \\
        &\mbox{and}&\quad  \delta^{-1}(n+1)^2e^{\frac{n-2}{2}}d'^{\frac{n+3}{2}}v^{\frac{n+3}{2}}r^{-2} = CL_f\delta^L \\
       & &\quad \Leftrightarrow \quad \delta =\gamma^{\frac{1}{L+1}} v^{\frac{n+3}{2(L+1)}}r^{\frac{-2}{L+1}}, \quad \mbox{where} \quad \gamma = \frac{(n+1)^2e^{\frac{n-2}{2}}d'^{\frac{n+3}{2}}}{CL_f}.
    \end{eqnarray*}
    Hence
    \begin{displaymath}
    bv^{-1} =\delta\quad \Rightarrow\quad v = \gamma^{\frac{-2}{n+2L+5}}b^{\frac{2(L+1)}{n+2L+5}}r^{\frac{4}{n+2L+5}}.
    \end{displaymath}
    Since $v$ is an integer, we make the following modification to the values of $\delta$ and $v$:
    \begin{equation}\label{choice of delta and v}
        v = \left\lfloor\gamma^{\frac{-2}{n+2L+5}} b^{\frac{2(L+1)}{n+2L+5}}r^{\frac{4}{n+2L+5}}\right \rfloor,\quad \delta =\gamma^{\frac{2}{n+2L+5}} b^{\frac{n+3}{n+2L+5}}r^{\frac{-4}{n+2L+5}}.
    \end{equation}
    Since $v$ is bounded by a fractional exponent of $r$, and $\delta$ is bounded by a negative fractional exponent of $r$, as $r$ tends to infinity, the conditions $r \geq 2n\max\{d,2l(v+1)\}$ and $\delta \in (0,1]$ are satisfied. 
    
    Substituting the values of $v$ and $\delta$ in \eqref{choice of delta and v} back to \eqref{approximation on bound} and using the inequality $\lfloor u \rfloor \leq u\ \forall u \in \bR$, we see that the first and third terms in \eqref{bound cpop} is $\mathcal{O}(mCL_f\delta^{L})$ and the second term is $\mathcal{O}(\delta^{-1}r^{-2})$. Whence, we obtain the following inequality:
    \begin{align}\label{bound cpop 2}
        &W(v,\delta,r) %\nonumber\\ =\;& 
        =\mathcal{O}\bigl(mCL_f \gamma^{\frac{2L}{n+2L+5}} b^{\frac{L(n+3)}{n+2L+5}}r^{\frac{-4L}{n+2L+5}}+m\gamma^{\frac{-2}{n+2L+5}} b^{\frac{-n-3}{n+2L+5}}r^{\frac{4}{n+2L+5}-2}\bigr)\nonumber\\
        &=\; \mathcal{O}\bigl(mCL_f \gamma^{\frac{2L}{n+2L+5}} \ln(r)^{\frac{L(n+3)}{n+2L+5}}r^{\frac{-4L}{n+2L+5}}+m\gamma^{\frac{-2}{n+2L+5}} \ln(r)^{\frac{-n-3}{n+2L+5}}r^{\frac{4}{n+2L+5}-2}\bigr).
    \end{align}
    The last equality is due to the fact that $b$ is proportional to $\ln(r)$ and $L \in (0,1]$. We next apply the inequality $\ln(x) < x^\alpha/\alpha\ \forall \alpha\in(0,1)$ and $x >1$ to obtain that $\ln(r) \leq (n+3)r^{\frac{1}{n+3}}$, which combines with \eqref{bound cpop 2} to give
    \begin{equation}\label{bound bpop 3}
        W(v,\delta,r) = \mathcal{O}\bigl(mnCL_f \gamma^{\frac{2L}{n+2L+5}} r^{\frac{-2L}{n+2L+5}}+m\gamma^{\frac{-2}{n+2L+5}}r^{\frac{-2L}{n+2L+5}}\bigr).
    \end{equation}
    Here, we used the inequality ${\frac{4}{n+2L+5}-2} \leq \frac{-2L}{n+2L+5}$
    and the fact that $\ln(r)^{\frac{-n-3}{n+2L+5}}\leq 1$ to simplify the second term in  \eqref{bound cpop 2}.
    By the definition of $\gamma$, we obtain that
    \begin{eqnarray*}
        CL_f\gamma^{\frac{2L}{n+2L+5}}&=& (CL_f)^{1-{\frac{2L}{n+2L+5}}}\bigl((n+1)^2e^{\frac{n-2}{2}}d'^{\frac{n+3}{2}}\bigr)^{\frac{2L}{n+2L+5}}\\
        &\leq & \max\{CL_f,1\}(n+1)^{\frac{4}{n+5}}e^{L}d'^{L}= \mathcal{O}\bigl( \max\{CL_f,1\}e^Ld'^L\bigr),\\
        \gamma^{\frac{-2}{n+2L+5}} &\leq& (CL_f)^{\frac{2}{n+2L+5}} \leq \max\{CL_f,1\}.
    \end{eqnarray*}
    Therefore, from \eqref{bound bpop 3}, we get
    \begin{align*}
        W(v,\delta,r) = \mathcal{O}\bigl(mn\max\{CL_f,1\}e^Ld'^L r^{\frac{-2L}{n+2L+5}}\bigr).
    \end{align*}
    In conclusion, there exists a positive constant $\rho$ independent of $n,m,l,L_f,C,L$ such that for any $r$ satisfying $r \geq 2n\max\{d,2l(v+1)\}$ and $\delta \in (0,1]$ with $v$ and $\delta$ defined by \eqref{choice of delta and v}, the following inequality holds:
    \begin{equation} \label{eq-rho}
        f_{\min}- \lb(f,\overline{\cQ}(\cX))_r \leq \rho 
        m n\max\{CL_f,1\}e^Ld'^L r^{\frac{-2L}{n+2L+5}}.
    \end{equation}
    Hence, the convergence rate of \eqref{proposed hierarchy for cpop} has a polynomial rate of $\mathcal{O}(r^{\frac{-2L}{n+2L+5}})$. 

    \medskip
    \noindent (ii) Next we consider the case where $f$ is positive over $\cX$.
    We aim to derive a bound on $r$ such that $f \in \widetilde{\cQ}(\cX)_{2r}$.  
    Recall that  $\overline{r} \geq 2n\max\{d,2l(v+1)\}$. For $\delta \in (0,1]$ and $v$ defined in \eqref{choice of delta and v},
    \eqref{bound bpop 3} is valid for any $r \geq \overline{r}$. 
    Furthermore, if $\lb(f,\overline{\cQ}(\cX))_r \geq 0$, then there exists $\sigma \in \overline{\cQ}(\cX)_{2r}$ such that 
    \begin{displaymath}
        f(x) = \sigma(x) + \lb(f,\overline{\cQ}(\cX))_r \in \overline{\cQ}(\cX)_{2r}.
    \end{displaymath}
    By \eqref{eq-rho},
    the inequality $\lb(f,\overline{\cQ}(\cX))_r \geq 0$ can be ensured for any $r \geq \overline{r}$ satisfying that 
    \begin{eqnarray*}
        &&\rho mn\max\{CL_f,1\}e^Ld'^L r^{\frac{-2L}{n+2L+5}} \leq f_{\min}=:E_f\\
        &\Leftrightarrow & \quad r \geq (\rho mn\max\{CL_f,1\})^{1+\frac{n+5}{2L}} (ed')^{\frac{n+2L+5}{2}}E_f^{-1-\frac{n+5}{2L}}.
    \end{eqnarray*}
    This completes the proof.
    \end{proof}

\begin{remark}
Theorem~\ref{thm: bound on degree} can be viewed as a generalization of Theorem~\ref{thm-Baldi} that provides bounds on the degrees of the SOS polynomials in the SOS representation associated with a matrix inequality. Moreover, this bound is valid for our proposed hierarchy \eqref{proposed hierarchy for cpop}, which is much simpler than the SOS representation proposed by Hol and Scherer.
\end{remark}
\section{Conclusion} 
    We have introduced a novel SOS hierarchy tailored for a polynomial optimization problem whose feasible set is defined by a matrix polynomial semidefinite inequality.
    The size of the moment matrices in the corresponding SDP relaxation is significantly smaller than the existing ones. We establish the convergence rate of the hierarchy through a penalty function approach. Our findings not only broaden the scope of scalar polynomial optimization to the wider matrix setting but also slightly improved the overall convergence rate. 
    %%Our future goal is to leverage analogous techniques to extend further advantageous attributes of the SOS hierarchy from scalar polynomial optimization to semidefinite polynomial optimization.
\section*{Acknowledgments}
We would like to express our sincere gratitude to the referees for their dedication and thorough review in providing numerous valuable feedback 
to improve the paper.

\appendix
\section{Construction of the univariate polynomials 
$T_k$ in Section~\ref{section: Penalty for binary} and related proofs}\label{appendix: univariate polynomial}
We first derive the properties of the function $\varphi(u)$ used in the proof
of Lemma~\ref{concatination conditon}. Then we provide the proofs of Proposition~\ref{prop: concatenation poly} and Lemma \ref{lemma: AC and BV} in Section~\ref{section: Penalty for binary}.

We refer the reader to \cite{rudin1953principles} for the definitions of absolute continuity, total variation, and functions of bounded variation.

\begin{lemma}\label{lem: Taylor expansion of varphi}
	The function $\varphi(u)$ is infinitely differentiable on $\left(-\infty,\frac{1}{4} \right)$, and its Taylor expansion at $u=0$ converges point-wise to $\varphi(u)$ on $(-1/4,1/4)$.
\end{lemma}
   \begin{proof}
      % For any $t <\frac{1}{2}$, $u = t(1-t) <\frac{1}{4}$ implies that 
      % \begin{displaymath}
      %     t = \frac{1}{2}-\sqrt{\frac{1}{4}-u} \quad \Rightarrow \quad \frac{t}{1-%2t}=\frac{2u}{1-4u + \sqrt{1-4u}}.
      % \end{displaymath}
%Thus, for any $u <\frac{1}{4}$, $\varphi(u) =\frac{2u}{1-4u + \sqrt{1-4u}}$. Moreover, 
Consider the function $\varphi$ defined on the domain $(-\infty, 1/4)$  by
\begin{displaymath}
   \varphi(u) =\; \frac{2u}{1-4u+\sqrt{1-4u}}=\; \frac{(1-\sqrt{1-4u})(1+\sqrt{1-4u})}{2\sqrt{1-4u}(1+\sqrt{1-4u})}=\; \frac{1}{2}(1-4u)^{-1/2}-\frac{1}{2}.
\end{displaymath}
By induction, we can prove that for any integer $k \geq 1$, the $k-$th derivative of $\varphi$ is 
\begin{equation}\label{derivative varphi in u}
\varphi^{(k)}(u)=\;\frac{1}{2}\prod_{i=0}^{k-1}\left(-\frac{1}{2}-i\right)(-4)^k(1-4u)^{-1/2-k}=\; \frac{1}{2}4^k\prod_{i=0}^{k-1}\left(i+\frac{1}{2}\right)(1-4u)^{-1/2-k}.
\end{equation}
Therefore, $\varphi$ is infinitely differentiable on $(-\infty,1/4)$, and its Taylor expansion at $u=0$ is given by 
\begin{equation}\label{Taylor expansion of varphi}
    \sum_{k=0}^{\infty}a_k u^k \quad \mbox{with}\quad a_0 =\varphi(0) =0, \;\; 
    a_k = \frac{\varphi^{(k)}(0)}{k!} = \frac{1}{2}\frac{4^k}{k!}\prod_{i=0}^{k-1}\left(i+\frac{1}{2}\right)
    \;\forall\; k\geq 1.
\end{equation}
The radius of convergence is calculated as 
\begin{displaymath}
    \left(\limsup_{k \to \infty} \sqrt[k]{\frac{1}{2}\frac{4^k}{k!}\prod_{i=0}^{k-1}\left(i+\frac{1}{2}\right)}\right)^{-1} \geq\; \left(\limsup_{k \to \infty} \sqrt[k]{\frac{1}{2}\frac{4^k}{k!}
    \frac{(2k-1)!!}{2^k}} \right)^{-1}= \frac{1}{4}.
\end{displaymath}
Whence, the Taylor expansion at $u=0$ converges point-wise to $\varphi(u)$ on $(-1/4,1/4)$.
\end{proof}

\begin{proof}[Proof of Proposition~\ref{prop: concatenation poly}.]
   For any positive integer $k$, consider the Taylor polynomial $T_k$ of $\varphi$ at $u=0$ given by \eqref{Taylor expansion of varphi} as 
   \begin{equation}\label{def: Taylor poly T_k}
       T_k(u)=\;\sum_{k=0}^{k}a_i u^i \quad \mbox{where}\quad 
       a_0 =0, \quad a_i = \frac{1}{2}\frac{4^i} {i!}\prod_{j=0}^{i-1}\left(j+\frac{1}{2}\right)\;>\; 0 \;\forall\; i\geq 1.
   \end{equation}
	Since $a_0=0$,  $T_k(u)$ has a factor of $u(t) = t(1-t)$. Whence, we obtain that $c_k(0)= 0$ and $c_k(1)=1$. Since $0 = t+(2t-1)\varphi(t(1-t))$ for all $t \in (-\infty,1/2)$, and for any $t \in [0,1/2)$, $u(t) = t(1-t) \in [0,1/4)$, we obtain that 
    for $c_k(t)=t+(2t-1)T_k(t(1-t))$,
    \begin{align*}
    c_k(t)&=\;c_k(t)-0 =\; (t+(2t-1)T_k(t(1-t))) - (t+(2t-1)\varphi(t(1-t))) \\
    &=\; (1-2t)\sum_{i=k+1}^{\infty}a_it^i(1-t)^{i} \quad \forall t \in [0,1/2).
    \end{align*}
    Therefore, $c_k^{(i)}(0) =0$ for all $i \in [k]$. We can see from \eqref{Taylor expansion of varphi}
  that $a_k$ is positive for all $k\geq 1$. Thus $T_k(u)$ monotonically increases to $\varphi(u)$ for any $u \in [0,1/4)$, and $T_k(u)$ is positive. Therefore, for any $t \in [0,1/2)$, we have  
    \begin{displaymath}
        0= t +(2t-1)\varphi(t(1-t)) \leq\;t+(2t-1)T_k(t(1-t))= c_k(t) \leq\; t \; \Rightarrow\;  0 \leq\; c_k(t) \leq\; \frac{1}{2}.
    \end{displaymath}
    Together with the fact that $c_k(t)+c_k(1-t) =1 \; \forall t \in [0,1]$, we obtain that $0 \leq c_k(t) \leq 1$ for all $t \in [0,1]$.
	  \end{proof}
   \begin{proof}[Proof of Lemma \ref{lemma: AC and BV}]
       Since $q_k(\lambda,N)$ is $k-$times continuously differentiable on the compact interval $[-1,1]$, then $q_k(\lambda,N)^{(k-1)}$ is absolutely continuous on $[-1,1]$ (see e.g., \cite{rudin1953principles}). Moreover, $q_k(\lambda,N)^{(k)}$ is a continuous piecewise polynomial given by 
       \begin{equation*}
	    q_k(\lambda,N)^{(k)}(t)=\begin{cases}
	  	0 &\mbox{if} \ t \in [-1,\lambda] \cup [0,1],\\
	  	N(c_k(t/\lambda))^{(k)}& \mbox{if}\ t \in [\lambda,0].
	  \end{cases}
	  \end{equation*}
      Thus, the total variation of $q_k(\lambda,N)^{(k)}$, denoted by $V_{[-1,1]}(q_k(\lambda,N)^{(k)})$, satisfies that
       \begin{multline}\label{total variation 1}
          V_{[-1,1]}(q_k(\lambda,N)^{(k)}) \leq\; \int_{\lambda}^0\left|N(c_k(t/\lambda))^{(k+1)}(t)\right|dt\\
          =\; N|\lambda|^{-k-1}\int_{\lambda}^0|c_k^{(k+1)}(t/\lambda)|dt=\;N|\lambda|^{-k}\int_{0}^1|c_k^{(k+1)}(t)|dt.
       \end{multline}
    We next bound $\int_{0}^1|c_k^{(k+1)}(t)|dt$. Recall that $u(t) = t(1-t)$.
    We apply the general Leibniz rule to compute the high order derivatives of $c_k$ as follows: for $k\geq 1$,
    \begin{eqnarray}\label{total variation 2}
        c_k^{(k+1)}(t) &=& t^{(k+1)}+[(2t-1)(T_k \circ u)(t) ]^{(k+1)}\nonumber\\
        &=&  (2t-1)(T_k \circ u)^{(k+1)}(t) + 2(k+1)(T_k \circ u)^{(k)}(t)\nonumber\\
        \Rightarrow \; \int_{0}^1|c_k^{(k+1)}(t)|dt\; &\leq&
        \int_{0}^1|(2t-1)(T_k \circ u)^{(k+1)}(t)|dt
        \nonumber \\
        && 
        \;+\;  2(k+1)\int_0^1|(T_k \circ u)^{(k)}(t)|dt.
    \end{eqnarray}
    We can bound the derivatives of $\varphi$ at $u=0$ 
    from \eqref{derivative varphi in u} as follows:
    \begin{equation}\label{bound on derivative of varphi}
        0 \leq \;\varphi^{(k)}(0)=\; \frac{1}{2}4^k\prod_{i=0}^{k-1}\left(i+\frac{1}{2}\right)\leq\; \frac{1}{2}4^kk!.
    \end{equation}
    We next calculate the high order derivatives of $T_k \circ u$ by Fa\`a di Bruno's formula (see, e.g., \cite{charalambides2018enumerative}) with Bell polynomials $B_{m,k}$ and the fact that $u(t) \in [0,1/4)$ for all $t \in [0,1]$. Furthermore, we note that the derivatives $u^{(i)}(t)$ vanishes for $i> 2$. Thus we only need to use the Bell polynomials $B_{k+1,k+1}$ and $B_{k+1,k}$ as follows: 
    \begin{eqnarray}\label{total variation 3}
        &(T_k \circ u)^{(k+1)}(t)&=\;T_k^{(k+1)}(u)B_{k+1,k+1}(u')+T_k^{(k)}(u)B_{k+1,k}(u')\nonumber\\
        &&=\;0+\varphi^{(k)}(0)\frac{k(k+1)}{2}(u')^{k-1}u''
        \quad \mbox{(by \eqref{def: Taylor poly T_k})}
        \nonumber\\
        \Rightarrow\;&\int_{0}^1|(2t-1)(T_k \circ u)^{(k+1)}(t)|dt \;
        &\leq\;\int_0^1\left|\varphi^{(k)}(0)\frac{k(k+1)}{2}(1-2t)^k(-2)\right|dt\nonumber\\
        &\;&\leq\;\frac{1}{2}4^k(k+1)!k\int_0^1|1-2t|^kdt \quad 
        \mbox{(by \eqref{bound on derivative of varphi})}
        \nonumber\\
        &&\leq\;\frac{1}{2}4^k(k+1)!k\frac{1}{k+1}=\;\frac{1}{2}4^kk!k.
    \end{eqnarray}
    Moreover,
    \begin{eqnarray}\label{total variation 4}
        (T_k \circ u)^{(k)}(t)&=&T_k^{(k)}(u)B_{k,k}(u')+T_k^{(k-1)}(u)B_{k,k-1}(u')\nonumber\\
        &=&\varphi^{(k)}(0)(u')^k+\Big(\varphi^{k-1}(0)+\frac{\varphi^{(k)}(0)}{k}u\Big)\frac{k(k-1)}{2}(u')^{k-2}u''
        \quad \mbox{(by \eqref{def: Taylor poly T_k})}
        \nonumber\\
        \Rightarrow\;\int_0^1|(T_k \circ u)^{(k)}(t)|dt\;&\leq& \left|\varphi^{(k)}(0)\right|\int_0^1\left|(u')^k\right|dt\nonumber\\
        &&+\;k(k-1)\left|\varphi^{k-1}(0)+\frac{\varphi^{(k)}(0)}{k}u\right|\int_0^1|u'|^{k-2}dt\nonumber\\
        \;&\leq& \frac{1}{2}4^kk!\int_0^1\left|1-2t\right|^kdt
        \quad \mbox{(by \eqref{bound on derivative of varphi}, and $u \in [0,1/4)$)}
        \nonumber\\
        &&+\;k(k-1)\left(\frac{1}{2}4^{k-1}(k-1)!+\frac{1}{2}\frac{4^kk!}{k}\frac{1}{4}\right)\int_0^1|1-2t|^{k-2}dt\nonumber\\
        &=& \frac{1}{2}4^kk!\frac{1}{k+1}+k(k-1)4^{k-1}(k-1)!\frac{1}{k-1}
         \;=\; \frac{1}{2}4^kk!\frac{1}{k+1}+4^{k-1}k!\nonumber\\
        &\leq & 2\cdot 4^{k-1}k!.
    \end{eqnarray}
    Substituting \eqref{total variation 3} and \eqref{total variation 4} into \eqref{total variation 2}, we get
    \begin{eqnarray}\label{total variation 5}
         &\int_{0}^1|c_k^{(k+1)}(t)|dt\; &\leq\; \frac{1}{2}4^kk!k+2(k+1)
         2 \cdot4^{k-1}k! \leq 3\cdot4^{k}k!k.
    \end{eqnarray}
    Substitute \eqref{total variation 5} into \eqref{total variation 1}, we obtain that 
        $V_{[-1,1]}(q_k(\lambda,N)^{(k)}) \leq\;3N|\lambda|^{-k}4^{k}k!k.$
   \end{proof}

%%%%%%%%%%%%%%%%%%%%%%%%%%%%%%%%%%%%   
\section{Inequalities} 
We derive the inequalities that we use
in the proofs in Sections~\ref{section: BPOP} and~\ref{section: Continuous}. 

\begin{lemma}\label{lem: choosing k}
    For any given positive integer $v$ and $\delta \in (0,1]$, 
    there exists a nonnegative integer $k<v$ such that 
    the following inequality holds:
    \begin{equation}\label{bound on the choice of k}
         \frac{4^kk!}{\delta^k(v-k)^k} \leq\; e^{2}a^{\delta v}, \quad \left(a= e^{-1/ (2e+1)}\approx 0.8561 \right).
    \end{equation}
\end{lemma}
\begin{proof}
 First, for any nonnegative integer $k$, we apply the Cauchy inequality to $k!$ to obtain
    \begin{displaymath}
        k\sqrt[k]{k!} \leq\; \sum_{i=1}^ki=\; \frac{k(k+1)}{2} \Rightarrow\; k! \leq\; \left(\frac{k+1}{2} \right)^k.
    \end{displaymath}
   Thus, we get
    \begin{equation}\label{eta function}
        \frac{4^kk!}{\delta^k(v-k)^k}\leq\; \frac{4^k(\frac{k+1}{2})^k}{\delta^k(v-k)^k}=\; \left(  \frac{2(k+1)}{\delta(v-k)}\right)^k =:\; \eta(k).
    \end{equation}
    Second, to bound $\eta(k)$, we approximate the minimizer of $\ln(\eta(k))$ 
    by considering the following derivative condition with respect to $k$:
    \begin{equation}\label{eq: derivative of ln}
       0= (\ln(\eta(k)))' =\; \left[k\ln\left(\frac{2(k+1)}{\delta(v-k)}\right)\right]'=\;\ln\left(\frac{2(k+1)}{\delta(v-k)}\right)+\frac{k}{k+1}+\frac{k}{v-k}. 
    \end{equation}
    Notice that there is no closed-form solution for the above equation. 
    To approximate the optimal $k$, we solve the following modified equation: 
    \begin{displaymath}
        \ln\left(\frac{2(k+1)}{\delta(v-k)}\right)+1=\;0 \quad \Leftrightarrow\quad \frac{2(k+1)}{\delta(v-k)}=\; \frac{1}{e} \quad \Leftrightarrow\quad k = \frac{\delta v -2e}{2e+\delta} < v.
    \end{displaymath}
    Thus we get that the approximately optimal $k$ is a small fraction of $v$ and it tends to $\infty$ as $v \to \infty$.
    Since $k$ should be a non-negative integer, we choose 
    \begin{equation} \label{eq-k}
        k=\; \max\left\{0,\; \left\lfloor \frac{\delta v -2e}{2e+\delta}\right\rfloor\right\}.
    \end{equation}
    Third, we substitute the above value of $k$ into $\dfrac{4^kk!}{\delta^k(v-k)^k}$. We consider two cases. (i) If  $\delta v -2e \geq 0$, then
    \begin{equation}\label{bound on base}
        k = \left\lfloor \frac{\delta v -2e}{2e+\delta}\right\rfloor \leq \frac{\delta v -2e}{2e+\delta} \Rightarrow\; \frac{2(k+1)}{\delta(v-k)} \leq\; e^{-1} < 1.
    \end{equation}
     Using the fact that $\delta \in (0,1]$, we also have
    \begin{equation}\label{bound on exponent}
        k \geq\; \frac{\delta v -2e}{2e+\delta}-1 =\; \frac{\delta v}{2e+\delta}-\frac{2e}{2e+\delta}-1 \geq\; \frac{\delta v}{2e+1}-2.
    \end{equation}
    Therefore, we combine \eqref{eta function}, \eqref{bound on base}, and \eqref{bound on exponent} to obtain that 
    \begin{equation*}\label{eq: inequality on k}
        \frac{4^kk!}{\delta^k(v-k)^k}\leq\;\left(  \frac{2(k+1)}{\delta(v-k)}\right)^k \leq\; e^{2-\delta v/ (2e+1)}=\; e^{2}a^{\delta v}.
    \end{equation*}
    (ii) If $\delta v -2e \leq 0$, then $k=0$, and 
    \begin{displaymath}
        \frac{4^kk!}{\delta^k(v-k)^k}=\;1, \quad \mbox{and}\quad e^{2-\delta v/ (2e+1)} \geq e^{2-2e/ (2e+1)} \geq e > 1.
    \end{displaymath}
    In both cases, the inequality \eqref{bound on the choice of k} holds
    for the value $k$ in \eqref{eq-k}.
\end{proof}

\begin{lemma}\label{lem: parameter C}
    Let $C_B(n,d)$ be the parameter stated in Theorem~\ref{unit ball}. Then for any positive integer $l,d$ and $v$, the following inequality holds:
    \begin{displaymath}
        C_B(n,\max\{d,2l(v+1)\})= \mathcal{O}\left((n+1)^2e^{\frac{n-2}{2}}d'^{\frac{n+3}{2}}v^{\frac{n+3}{2}}\right),
    \end{displaymath}
    where $d' = \max\{ \lceil d/2\rceil, l\}$.
\end{lemma}
\begin{proof}
    The explicit expression of $C_B(n,d)$ is stated in the paper \cite{slot2111sum}. We summarize the main results as follows: for any positive integers $n$ and $d$,
    \begin{displaymath}
        C_B(n,d) = 2(n+1)^2d^2\gamma(B^n)_d, \quad \mbox{and} \quad \gamma(B^n)^2_d=\max_{0 \leq k \leq d}\left(1+\frac{2k}{n-1}\right)\cdot\binom{k+n-2}{n-2}.
    \end{displaymath}
    For any $0 \leq k \leq d$,  using the inequality $\binom{n}{k} \leq \left(\frac{en}{k} \right)^k$, we have that
    \begin{eqnarray*}
        &\binom{k+n-2}{n-2}& \leq \left(\frac{e(k+n-2)}{n-2} \right)^{n-2} \leq \left(e + \frac{ed}{n-2} \right)^{n-2}= e^{n-2}\left(1+\frac{d}{n-2} \right)^{n-2}\\
        &\Rightarrow& \gamma(B^n)_d \leq e^{\frac{n-2}{2}}\left(1+\frac{2d}{n-1} \right)^{1/2}\left(1+\frac{d}{n-2} \right)^{\frac{n-2}{2}}.
    \end{eqnarray*}
    Note that $2d'(v+1) > d$. We can
    evaluate
    \begin{align*}
        &C_B(n,\max\{d,2l(v+1)\}) \leq C_B(n,2d'(v+1))\\
        \leq &\ 2(n+1)^2(2d'(v+1))^2e^{\frac{n-2}{2}}\left(1+\frac{4d'(v+1)}{n-1} \right)^{1/2}\left(1+\frac{2d'(v+1)}{n-2} \right)^{\frac{n-2}{2}}\\
        \leq &\ 2(n+1)^2(2d'(v+1))^2e^{\frac{n-2}{2}}\left(1+\frac{4d'(v+1)}{n-1} \right)^{1/2}(d'v)^{\frac{n-2}{2}}\left(1+\frac{4}{n-2} \right)^{\frac{n-2}{2}}\\
        =&\ \mathcal{O}\left((n+1)^2e^{\frac{n-2}{2}}d'^{\frac{n+3}{2}}v^{\frac{n+3}{2}}\right), \quad \text{since} \lim_{n \to \infty}\left(1+\frac{4}{n-2} \right)^{\frac{n-2}{2}}= e^2.
    \end{align*}
\end{proof}

\bibliographystyle{siamplain}
\bibliography{references}
\end{document}

%% file: ex_shared.tex
\usepackage{verbatim}
\usepackage{amsfonts}
\usepackage{graphicx}
\usepackage{epstopdf}
\usepackage{subcaption}
\usepackage{amssymb}
\usepackage{algorithmic}
\ifpdf
  \DeclareGraphicsExtensions{.eps,.pdf,.png,.jpg}
\else
  \DeclareGraphicsExtensions{.eps}
\fi

\newsiamremark{remark}{Remark}
\newsiamremark{hypothesis}{Hypothesis}
\crefname{hypothesis}{Hypothesis}{Hypotheses}
\newsiamthm{claim}{Claim}
\newtheorem{assume}{Assumption}
\newsiamremark{example}{Example}

\allowdisplaybreaks

\title{Convergence rates of SoS hierarchies for polynomial semidefinite programs}
\author{Hoang Anh Tran, \ 
  Kim-Chuan Toh}

\usepackage{amsopn}
\DeclareMathOperator{\diag}{diag}

\newcommand{\bR}{\mathbb{R}}
\newcommand{\bN}{\mathbb{N}}
\newcommand{\bS}{\mathbb{S}}
\newcommand{\cQ}{\mathcal{Q}}
\newcommand{\cX}{\mathcal{X}}
\newcommand{\tr}{\mathrm{tr}}
\newcommand{\bB}{\mathbb{B}}
\newcommand{\M}{\mathbf{M}}
\newcommand{\cC}{\mathcal{C}}
\newcommand{\cP}{\mathcal{P}}
\newcommand{\cH}{\mathcal{H}}